\documentclass[12pt]{article}
\usepackage{amsmath,amssymb,amsthm}
\usepackage{color}
\usepackage{comment}
\numberwithin{equation}{section}
\newtheorem{thm}{Theorem}[section]
\newtheorem{prop}[thm]{Proposition}
\newtheorem{cor}[thm]{Corollary}
\newtheorem{exam}[thm]{Example}
\newtheorem{rem}[thm]{Remark}
\newtheorem{lem}[thm]{Lemma}
\newtheorem{defn}[thm]{Definition}
\newtheorem{assum}[thm]{Assumption}
\newcommand{\pf}{{\em Proof.}}

\def\address#1#2{\begingroup
\noindent\parbox[t]{16cm}{%
\small{\scshape\ignorespaces#1}\par\vskip1ex
\noindent\small{\itshape E-mail address}%
\/: #2\par\vskip4ex}\hfill%
\endgroup}%

\begin{document}
\title{Maximal displacement  and population growth \\ 
for branching Brownian motions}
\author{Yuichi Shiozawa
\footnote{Supported in part by JSPS KAKENHI No.\ JP17K05299.}}
\maketitle 

\begin{abstract}
We study the maximal displacement and related population 
for a branching Brownian motion in Euclidean space 
in terms of the principal eigenvalue of an associated Schr\"odinger type operator.  
We first determine their growth rates on the survival event.  
We then establish the upper deviation for the maximal displacement 
under the possibility of extinction. 
Under the non-extinction condition, 
we further discuss the decay rate of the upper deviation probability  
and the population growth at the critical phase. 

\end{abstract}

\maketitle
\section{Introduction}

We are concerned with the population growth rate related to the {\it maximal displacement} 
for a {\it spatially inhomogeneous} branching Brownian motion in Euclidean space ${\mathbb R}^d$.
We proved in \cite{S18} that under the non-extinction condition, 
this rate is given  
in terms of the principal eigenvalue of an associated Schr\"odinger type operator. 
This result implies the existence of the phase transition for the growth rate. 
As its corollary, we determined the linear growth rate of the maximal displacement. 
We further established the upper deviation for the maximal displacement. 
In this paper, we first remove the non-extinction condition in \cite{S18} 
(Theorem \ref{thm:main}, Corollary \ref{cor:main} and Theorem \ref{thm:upper}). 
We next discuss the decay rate of the tail probability of the maximal displacement 
as a refinement of the upper deviation under the non-extinction condition (Theorem \ref{thm:decay-rate}). 
We finally prove that for $d\geq 3$, 
the population growth rate as mentioned before is polynomial 
at the critical phase 
under the same condition (Theorem \ref{thm-critical}).

The maximal displacement is one of the important research objects for branching Brownian motions 
because it reflects quantitatively the interplay 
between the randomness of branching and that of particle motions. 
The distribution of the maximal displacement is also related 
to the so called Fisher-Kolmogorov-Petrovskii-Piskunov equation 
(see, e.g., \cite{CR88, E84, Mc75, Mc76, S18} and references therein).  
We would like to mention some of the results, 
which are related to the problems in this paper,  
for a one dimensional binary branching Brownian motion 
such that the splitting time is  exponentially distributed with rate $c>0$. 
As is well known, the maximal displacement $R_t$ at time $t$ satisfies 
the law of large numbers  $R_t/t\rightarrow \sqrt{2c} \ (t\rightarrow\infty)$ a.s.\ 
(see, e.g., Bramson \cite{B78} and Roberts \cite{R13} for more detailed properties). 
Chauvin and Rouault \cite{CR88} determined 
the decay rate of the probability of the upper deviation type 
for the maximal displacement.  
Biggins \cite{Bi95,Bi96} further obtained the growth rate 
of the population right to the point $\delta t$ at time $t$, 
where $\delta$ is a positive constant such that $\delta\ne \sqrt{2c}$. 
We note that the law of large numbers for the maximal displacement is valid  
also for $d\geq 2$ and the offspring distribution is more general 
so that extinction may occur; however, that distribution is assumed to be spatially independent
(see, e.g., \cite{K05}, \cite{M15}, \cite{OCE17}). 
Biggins \cite{Bi95,Bi96} also mentioned that his result is valid 
under a setting similar to that as above. 

Our interest here is how the spatial inhomogeneity of the branching structure 
affects the behavior of the population growth related to the maximal displacement. 
By the spatial inhomogeneity,  
we mean that the distributions of the splitting time and offspring 
depend on the trajectory of each particle and branching site, respectively 
(see Subsection \ref{subsect-model} below for details). 
As for the population size,  
the long time behavior is characterized in terms of the principal eigenvalue of a Schr\"odinger type 
operator associated with the branching structure 
(see, e.g., \cite{CRY17, CS07, EHK10, KM13, W67}). 
This characterization also applies to the maximal displacement.  
In fact, when $d=1$ and non-extinction occurs,
Erickson \cite{E84} proved that 
even if the branching intensity is small at infinity, 
the maximal displacement grows linearly and its rate is determined by 
the same principal eigenvalue as mentioned before. 
This result is valid also for $d\geq 2$ if the branching intensity is spherically symmetric. 
We can further obtain the exponential growth rate 
of the population outside balls with time dependent radius for $d\geq 1$ 
under the setting similar to that in \cite{E84}
(see \cite{BH14}, \cite{KM13}, \cite{S18}). 
In particular, we can allow the spherical asymmetry of the branching intensity. 
This result is regarded as a spatially inhomogeneous counterpart of Biggins \cite{Bi95, Bi96}. 
We note that the results of \cite{E84} and \cite{BH14} 
are also extended by Lalley and Sellke \cite{LS88} and Bocharov and Wang \cite{BW18+}, 
respectively, to the model in which 
the branching intensity is inhomogeneous and not small at infinity.

In connection with the extinction problem, 
it is natural to allow the possibility of extinction for the spatially inhomogeneous model. 
More precisely, we would like to see 
the behavior of the maximal displacement under the {\it survival event}. 
Our results (Theorem \ref{thm:main}, Corollary \ref{cor:main} and Theorem \ref{thm:upper}) 
say that the previous results in \cite{S18} remain true, 
and the effect of the possibility of extinction appears 
in the principal eigenvalue of the Schr\"odinger type operator as mentioned before. 
Our approach is similar to that of \cite{S18}, 
which is an extension of \cite{BH14} 
to the multidimensional branching Brownian motions with singular branching intensity; 
however, we overcome several difficulties arising  
from the fact that the total population is not increasing in time 
(see comments just after Lemmas \ref{lem:upper-growth} and \ref{lem:lower-growth}). 
We also reveal the long time behavior of the expected 
Feynman-Kac functional associated with a {\it signed} measure 
(see \eqref{eq-br-fk} and comment just after Theorem \ref{thm-fk-g} below). 

Corollary \ref{cor:main} is partially regarded  
as a continuous (time-)space counterpart of 
Carmona and Hu \cite{CH14} and Bulinskaya \cite{B18}.  
They studied the growth rate of the maximal displacement 
for a branching random walk on the integer lattice 
such that each particle moves as a general irreducible (non-symmetric) random walk 
and branching occurs only on finite points. 
They also allow the possibility of extinction. 
As for our model, even though we assume that each particle moves as a Brownian motion, 
branching can occur on a non-compact set. 

Our refinement on the upper deviation type probability 
of the maximal displacement (Theorem \ref{thm:decay-rate}) 
is regarded as a spatially inhomogeneous counterpart of Chauvin and Rouault \cite{CR88}. 
In particular, we determine the exponential decay rate of this probability 
more precisely than \cite{S18}, and bound the polynomial order.    
Our argument is also similar to that of \cite{CR88}. 
For the lower bound of the probability especially, 
we utilize its Feynman-Kac expression 
originating from McKean \cite{Mc75, Mc76} (see \eqref{eq:fk-1} and \eqref{eq:fk-ex}).    
Here we impose the non-extinction condition on the branching structure  
because of the inequality \eqref{eq:side} below. 
We do not know if this condition can be dropped.

Theorem \ref{thm-critical} provides an information about 
the long time behavior of the population around the forefront. 
In particular, we see that for $d\geq 3$,  
such population grows polynomially with dimension dependent  growth rate. 
Our approach for Theorem \ref{thm-critical} is a refinement of 
that applied to the non-critical case in \cite{S18}.  
To derive the polynomial growth, 
we make use of the long time behavior of the Feynman-Kac functional 
associated with a {\it positive} measure (see \eqref{eq:fk-asymp-f} below). 
This also imposes the non-extinction condition on the branching structure. 
To the best of the author's knowledge,   
there are no references on the population growth around the forefront.

The rest of this paper is organized as follows. 
In Section \ref{sect:prelim}, 
we first introduce the Kato class measure and Feynman-Kac semigroups. 
We then introduce the model of branching Brownian motions. 
In Section \ref{sect:result}, we present our results in this paper 
and their applications to some concrete models. 
In Section \ref{sect:growth}, we derive the exponential growth rate  
of the expectation of the Feynman-Kac functional 
associated with a signed measure.  
The subsequent sections are devoted to the proofs of the results 
presented in Section \ref{sect:result}.
In Appendix \ref{appendix:decay},
we show a convergence result for the expectation of the Feynman-Kac functional
associated with a signed measure (see \eqref{eq-fkf-asymp} below).  
We follow the argument of Carmona \cite{C89} and Takeda \cite[Theorem 5.2]{T08}. 
In Appendix \ref{appendix:positive}, 
we discuss the relation between the regular growth and survival 
in order for the validity of the consequence of 
Theorem \ref{thm:main} and Corollary \ref{cor:main} 
on the survival event (see Remark \ref{rem:survival}). 
In Appendix \ref{appendix:evaluate}, 
we give a part on the elementary calculation in Section \ref{sect:growth}.

Throughout this paper, the letters $c$ and $C$ (with subscript) 
denote finite positive constants which may vary from place to place. 
For positive functions $f(t)$ and $g(t)$ on $(0,\infty)$, we write 
$f(t)\asymp g(t) \ (t\rightarrow\infty)$ if there exist positive constants $T$, $c_1$ and $c_2$ such that 
$c_1g(t)\leq f(t)\leq c_2g(t)$ for all $t\geq T$.  
We also write $f(t)\lesssim  g(t) \ (t\rightarrow \infty)$ if there exist positive constants $T$ and  $c$ such that 
$f(t)\leq cg(t)$ for all $t\geq T$.

\section{Preliminaries}\label{sect:prelim}
\subsection{Kato class measures and Feynman-Kac semigroups}\label{subsection:fk}
Let ${\bf M}=(\Omega, {\cal F}, \{{\cal F}_t\}_{t\geq 0}, \{B_t\}_{t\geq 0}, \{P_x\}_{x\in {\mathbb R}^d}, 
\{\theta_t\}_{t\geq 0})$ 
be the Brownian motion on ${\mathbb R}^d$, 
where $\{{\cal F}_t\}_{t\geq 0}$ is the minimal admissible filtration 
and $\{\theta_t\}_{t\geq 0}$ is the time shift operator of the paths 
such that $B_s\circ\theta_t=B_{s+t}$ identically for $s,t\geq 0$. 
Let 
$$p_t(x,y)=\frac{1}{(2\pi t)^{d/2}}\exp\left(-\frac{|x-y|^2}{2t}\right) \quad 
(x,y\in {\mathbb R}^d, t>0).$$
Then $p_t(x,y)$ is  the density of the transition function of ${\bf M}$, that is, 
$$P_x(B_t\in A)=\int_A p_t(x,y)\,{\rm d}y, \quad A\in {\cal B}({\mathbb R}^d).$$
Here $ {\cal B}({\mathbb R}^d)$ is the totality  of Borel subsets of ${\mathbb R}^d$. 
For $\alpha>0$, 
let $G_{\alpha}(x,y)$ be the $\alpha$-resolvent density of ${\bf M}$: 
$$G_{\alpha}(x,y)=\int_0^{\infty}e^{-\alpha t}p_t(x,y)\,{\rm d}t.$$
Then for any $\alpha>0$, 
\begin{equation}\label{eq:resolvent}
G_{\alpha}(x,y)\sim c\frac{e^{-\sqrt{2\alpha}|x-y|}}{|x-y|^{(d-1)/2}} \quad (|x-y|\rightarrow\infty)
\end{equation}
(see, e.g., \cite[(2.1)]{S18}).
For $d\geq 3$, we denote by $G(x,y)$ the Green function of ${\bf M}$:
$$G(x,y)=\int_0^{\infty}p_t(x,y)\,{\rm d}t=\frac{\Gamma(d/2-1)}{2\pi^{d/2}}\frac{1}{|x-y|^{d-2}}.$$
We also define $G_0(x,y):=G(x,y)$. 

According to \cite{ABM91, C02, T08},
we first introduce two classes of measures:

\begin{defn}\label{def-kato}
\begin{enumerate}
\item[{\rm (1)}] Let $\mu$ be a positive Radon measure on ${\mathbb R}^d$. 
Then $\mu$ belongs to the Kato class {\rm (}$\mu\in {\cal K}$ in notation{\rm )} 
if one of the following holds:
\begin{enumerate}
\item[{\rm (i)}] $d=1$ and 
$$\sup_{x\in {\mathbb R}}\int_{|x-y|\leq 1}\mu({\rm d}y)<\infty{\rm ;}$$
\item[{\rm (ii)}] $d=2$ and 
$$\lim_{R\rightarrow +0}\sup_{x\in {\mathbb R}^2}\int_{|x-y|\leq R}
\log\left(\frac{1}{|x-y|}\right)\,\mu({\rm d}y)=0{\rm ;}$$
\item[\rm (iii)] $d\geq 3$ and 
$$\lim_{R\rightarrow +0}\sup_{x\in {\mathbb R}^d}\int_{|x-y|\leq R}
G(x,y)\mu({\rm d}y)=0.$$
\end{enumerate}

\item[{\rm (2)}] For $\beta>0$, $\mu\in {\cal K}$ is $\beta$-Green tight 
{\rm (}$\mu\in {\cal K}_{\infty}(\beta)$ in notation{\rm)} if 
$$\lim_{R\rightarrow\infty}\sup_{x\in {\mathbb R}^d}\int_{|y|\geq R}G_{\beta}(x,y)\,\mu({\rm d}y)=0.$$
When $d\geq 3$, $\mu\in {\cal K}$ is Green tight if the equality above is valid 
for $\beta=0$.  
\end{enumerate}
\end{defn}

We know by  \cite{T08} and \cite[Corollary 4.2 and Lemma 4.2]{TTT17} that 
${\cal K}_{\infty}(\beta) \ (\beta>0)$ is independent of $\beta$ 
and ${\cal K}_{\infty}(0)\subsetneq {\cal K}_{\infty}(1)$. 
Define
$$
{\cal K}_{\infty}=
\begin{cases}{\cal K}_{\infty}(1) & (d=1,2), \\
{\cal K}_{\infty}(0) & (d\geq 3).
\end{cases}$$
If $\mu$ is a Kato class measure with compact support in ${\mathbb R}^d$, 
then $\mu\in {\cal K}_{\infty}$ by definition. 
For examples of measures in ${\cal K}_{\infty}$, 
see \cite[Subsection 2.1]{S18} and references therein.

We next introduce the notion of  positive continuous additive functionals.
Let $A=\{A_t\}_{t\geq 0}$ be a $[0,\infty]$-valued stochastic process on $(\Omega,{\cal F})$. 
We say that $A$ is a {\it positive continuous additive functional {\rm (}in the strict sense{\rm )}}  
(PCAF in short) of ${\mathbf M}$ if 
\begin{enumerate}
\item $A_t$ is ${\cal F}_t$-measurable for any $t\geq 0$;
\item There exists an event $\Lambda\in {\cal F}_{\infty}$, which satisfies 
$P_x(\Lambda)=1$ for any $x\in {\mathbb R}^d$ and 
$\theta_t\Lambda\subset \Lambda$ for any $t>0$, 
such that for any $\omega\in \Lambda$, 
\begin{itemize}
\item $A_0(\omega)=0$;
\item $A_t(\omega)$ is finite and continuous in $t\in [0,\infty)$;
\item $A_{t+s}(\omega)=A_t(\omega)+A_s(\theta_t\omega)$ for any $s,t\geq 0$
\end{itemize}
\end{enumerate}
(see, e.g., \cite{ABM91} and \cite[p.401]{FOT11}).
For each $\mu\in {\cal K}$, there exists 
a unique PCAF ($A^{\mu}$ in notation) 
such that for any nonnegative Borel function $f$, 
$$\lim_{t\rightarrow 0}\frac{1}{t}\int_{{\mathbb R}^d}E_x\left[\int_0^t f(B_s)\,{\rm d}A_s^{\mu}\right]\,{\rm d}x
=\int_{{\mathbb R}^d}f(x)\,\mu({\rm d}x)$$
(\cite[Proposition 3.8]{ABM91} and  \cite[Theorems 5.1.3 and  5.1.7]{FOT11}). 
We note that if $d\geq 3$, 
then by \cite[Proposition 2.2]{C02}, any measure $\mu\in {\cal K}_{\infty}$ is {\it Green-bounded}:
\begin{equation}\label{eq-green-bd}
\sup_{x\in {\mathbb R}^d}
E_x\left[A_{\infty}^{\mu}\right]=\sup_{x\in {\mathbb R}^d} \int_{{\mathbb R}^d}G(x,y)\,\mu({\rm d}y)<\infty.
\end{equation}

Let $\mu$ be a signed measure on ${\mathbb R}^d$ such that $\mu=\mu^{+}-\mu^{-}$ 
for some $\mu^+, \mu^-\in {\cal K}$ and define  $A_t^{\mu}=A_t^{\mu^+}-A_t^{\mu^-}$.
Then the multiplicative functional $e^{A_t^{\mu}}$ is called the {\it Feynman-Kac functional}. 
Using it, we  define the {\it Feynman-Kac semigroup} $\{p_t^{\mu}\}_{t\geq 0}$ by 
$$p_t^{\mu}f(x):=E_x\left[e^{A_t^{\mu}}f(B_t)\right], \quad f\in L^2({\mathbb R}^d)\cap {\cal B}_b({\mathbb R}^d),$$
where ${\cal B}_b({\mathbb R}^d)$ stands for the set of bounded Borel functions on ${\mathbb R}^d$.
Then $\{p_t^{\mu}\}_{t\geq 0}$ forms  
a strongly continuous symmetric semigroup on $L^2({\mathbb R}^d)$ 
such that its $L^2$-generator is formally expressed  as 
the Schr\"odinger type operator ${\cal H}^{\mu}:=-\Delta/2-\mu$.  
We can further extend $\{p_t^{\mu}\}_{t\geq 0}$ to $L^p({\mathbb R}^d)$ 
for any $p\in [1,\infty]$ (\cite[Theorem 6.1 (i)]{ABM91}). 
For simplicity, we use the same notation for the extended semigroup. 
Let $\|\cdot\|_{p,q}$ denote the operator norm from $L^p({\mathbb R}^d)$ to $L^q({\mathbb R}^d)$. 
We then have

\begin{thm}\label{thm-abm} {\rm (\cite[Theorems 6.1 (iii) and 7.1 (ii)]{ABM91})} 
Let $\mu$ be a signed measure on ${\mathbb R}^d$ such that 
$\mu=\mu^+-\mu^-$ for some $\mu^+, \mu^- \in {\cal K}$.
\begin{enumerate}
\item For any $t>0$, $\|p_t^{\mu}\|_{p,q}<\infty$ for any $1\leq p\leq q\leq \infty$.
\item For any $f\in {\cal B}_b({\mathbb R}^d)$ and $t>0$, $p_t^{\mu}f$ is a bounded continuous function on ${\mathbb R}^d$.
\end{enumerate}
\end{thm}

Assume that $\mu=\mu^+-\mu^-$ for some $\mu^+, \mu^- \in {\cal K}_{\infty}(1)$.   
Define 
$$\lambda({\mu}):=\inf\left\{\frac{1}{2}\int_{{\mathbb R}^d}|\nabla u|^2\,{\rm d}x-\int_{{\mathbb R}^d}u^2\,{\rm d}\mu 
\mid u\in C_0^{\infty}({\mathbb R}^d), \int_{{\mathbb R}^d}u^2\,{\rm d}x=1\right\},$$
where $C_0^{\infty}({\mathbb R}^d)$ stands for the set of smooth functions on ${\mathbb R}^d$ with compact support. 
Then $\lambda(\mu)$ is the bottom of the $L^2$-spectrum of ${\cal H}^{\mu}$. 
In particular,  
if $\lambda(\mu)<0$, then $\lambda(\mu)$ is the principal eigenvalue of ${\cal H}^{\mu}$ 
(see \cite[Lemma 4.3]{T03} or \cite[Theorem 2.8]{T08}) 
and the corresponding eigenfunction $h$ 
has a bounded, continuous and strictly positive version by  Theorem \ref{thm-abm}
(see, e.g., \cite[Section 4]{T08}).

In what follows, we assume that $\lambda:=\lambda(\mu)<0$ and the eigenfunction $h$ is 
bounded, continuous and strictly positive on ${\mathbb R}^d$
such that $\int_{{\mathbb R}^d}h(x)^2\,{\rm d}x=1$. 
Then by the proof of \cite[Theorem 5.2]{T08} (see also Subsection \ref{appendix:decay}), 
we have for any $f\in {\cal B}_b({\mathbb R}^d)$,
\begin{equation}\label{eq-fkf-asymp}
\lim_{t\rightarrow\infty}e^{\lambda t}E_x\left[e^{A_t^{\mu}}f(B_t)\right]
=h(x)\int_{{\mathbb R}^d}f(y)h(y)\,{\rm d}y \quad (x\in {\mathbb R}^d).
\end{equation}

\subsection{Branching Brownian motions}\label{subsect-model}
In this subsection, we introduce the model of branching Brownian motions 
by following \cite{INW68-1,INW68-2, INW69}. 
For $x\in {\mathbb R}^d$, let $\{p_n(x)\}_{n\geq 0}$ be a sequence such that 
$$0\leq p_n(x)\leq 1  \quad (n\geq 0) \quad \text{and} 
\quad \sum_{n=0}^{\infty}p_n(x)=1.$$
Let $\tau$ be the nonnegative random variable defined on $(\Omega, {\cal F}, P_x)$, 
which is independent of the Brownian motion, of exponential distribution 
with rate $1$; $P(\tau>t)=e^{-t}$ for any $t>0$. 
Let $\mu$ be a Kato class measure on ${\mathbb R}^d$.
We define 
$$Z:=\inf\left\{t>0 \mid A_t^{\mu}\geq \tau\right\}$$
so that 
$$P_x(Z>t \mid {\cal F}_{\infty})=e^{-A_t^{\mu}}.$$

We can describe the branching Brownian motion as follows: 
a Brownian particle $\{B_t\}_{t\geq 0}$ starts at $x\in {\mathbb R}^d$ 
according to the law $P_x$. 
At time $Z$, this particle splits into $n$ particles with probability $p_n(B_Z)$. 
These particles then start at $B_Z$ independently according to the law $P_{B_Z}$, 
and each of them continues the same procedure.

Let $({\mathbb R}^d)^{(0)}=\{\Delta\}$ and $({\mathbb R}^d)^{(1)}={\mathbb R}^d$. 
For $n\geq 2$, we define the equivalent relation $\sim$ 
on $({\mathbb R}^d)^n=\underbrace{{\mathbb R}^d\times \dots\times {\mathbb R}^d}_{n}$ as follows:
for ${\mathbf x}^n=(x^1,\dots, x^n)$ and ${\mathbf y}^n=(y^1,\dots, y^n)\in ({\mathbb R}^d)^n$, 
we write ${\mathbf x}\sim {\mathbf y}$ if there exists a permutation $\sigma$ on $\{1,2,\dots, n\}$ 
such that $y^i=x^{\sigma(i)}$ for any $i\in \{1,2,\dots, n\}$. 
If we define  $({\mathbb R}^d)^{(n)}=({\mathbb R}^d)^n/\sim$ for $n\geq 2$
and ${\mathbf X}=\cup_{n=0}^{\infty}({\mathbb R}^d)^{(n)}$, 
then $n$ points in ${\mathbb R}^d$ determine a point in $({\mathbb R}^d)^{(n)}$.
Hence we can define 
a branching Brownian motion 
$\overline{\mathbf M}=(\{{\mathbf B}_t\}_{t\geq 0}, \{{\mathbf P}_{\mathbf x}\}_{{\mathbf x}\in {\mathbf X}})$ 
on ${\mathbf X}$ with branching rate $\mu$ and branching mechanism $\{p_n(x)\}_{n\geq 0}$.  

Let $T$ be the first splitting time of $\overline{\mathbf M}$ given by 
\begin{equation}\label{eq-split}
{\mathbf P}_x(T>t \mid \sigma(B))=P_x(Z>t \mid {\cal F}_{\infty})=e^{-A_t^{\mu}} \quad (t>0).
\end{equation}
By definition, the first splitting time becomes small if the particle 
moves on the support of $\mu$ often. 
Let 
$$Q(x):=\sum_{n=0}^{\infty}np_n(x)$$
be the expected offspring number at branching site $x\in {\mathbb R}^d$. 
Denote by $Z_t$ the total number of particles at time $t$, that is, 
$$Z_t=n \quad 
\text{if ${\mathbf B}_t=({\mathbf B}_t^{1},\dots, {\mathbf B}_t^n)
\in {({\mathbb R}^d)}^{(n)}$.}$$
Let 
$$e_0=\inf\left\{t>0 \mid Z_t=0\right\}$$
be the extinction time of $\overline{\mathbf M}$ 
and $u_e(x)={\mathbf P}_x(e_0<\infty)$. 
We say that $\overline{\mathbf M}$ becomes extinct if $u_e\equiv 1$.

We define for $f\in {\cal B}_b({\mathbb R}^d)$, 
$$Z_t(f):=\sum_{k=1}^{Z_t}f({\mathbf B}_t^k).$$
For $A\in {\cal B}({\mathbb R}^d)$, 
let $Z_t(A):=Z_t({\bf 1}_A)$ be the number of particles on the set $A$ at time $t$.
If the measure 
$$
\nu({\rm d}x):=(Q(x)-1)\mu({\rm d}x)
$$
is written as $\nu=\nu^{+}-\nu^{-}$ for some $\nu^{+}, \nu^{-}\in {\cal K}$, 
then by the same way as in \cite[Lemma 3.3]{S08}, 
we have 
\begin{equation}\label{eq-br-fk}
{\mathbf E}_x\left[Z_t(f)\right]
=E_x\left[e^{A_t^{\nu}}f(B_t)\right].
\end{equation}

Assume that $\nu^+, \nu^-\in {\cal K}_{\infty}(1)$ and $\lambda:=\lambda(\nu)<0$. 
Let $h$ be the eigenfunction of ${\cal H}^{(Q-1)\mu}$ corresponding to $\lambda$ and 
$$M_t:=e^{\lambda t}Z_t(h) \quad (t\geq 0).$$
Since $M_t$ is a nonnegative ${\mathbf P}_x$-martingale,
the limit $M_{\infty}:=\lim_{t\rightarrow\infty}M_t\in [0,\infty)$ exists 
${\mathbf P}_x$-a.s. 
Furthermore, by \cite[Theorem 3.7]{CS07}, there exists an event of ${\mathbf P}_x$-full probability measure 
such that we have on this event,   
\begin{equation}\label{eq-limit-thm}
\lim_{t\rightarrow\infty}e^{\lambda t}Z_t(A)=M_{\infty}\int_A h(y)\,{\rm d}y
\end{equation}
for any $A\in {\cal B}({\mathbb R}^d)$ 
such that its boundary has zero Lebesgue measure.

\section{Results}\label{sect:result}
In this section, we state the results in this paper.
Let $\overline{\mathbf M}=(\{{\mathbf B}_t\}_{t\geq 0}, \{P_{\mathbf x}\}_{{\mathbf x}\in {\mathbf X}})$ 
be a branching Brownian motion on ${\mathbf X}$ 
with branching rate $\mu$ and branching mechanism $\{p_n(x)\}_{n\geq 0}$.  
We impose the next assumption on the branching rate and mechanism.
\begin{assum}\rm Let $\nu({\rm d}x)=(Q(x)-1)\mu({\rm d}x)$ and $\lambda=\lambda(\nu)$.
\begin{enumerate}\label{assum:1}
\item[(i)] $\lambda<0$.
\item[(ii)] The measure $\nu$ is written as 
$\nu=\nu^{+}-\nu^{-}$ for some $\nu^{+},\nu^{-}$ such that  
for any $\beta>0$, $\nu_{\beta}^+({\rm d}x):=e^{\beta |x|}\,\nu^+({\rm d}x)$ 
and $\nu_{\beta}^-({\rm d}x):=e^{\beta |x|}\,\nu^{-}({\rm d}x)$ belong to ${\cal K}_{\infty}(1)$. 
\end{enumerate}
\end{assum}
The condition (ii) says that the measure $\nu$ is small enough at infinity. 
In particular, this condition implies that $\nu^{\pm}({\mathbb R}^d)<\infty$ 
because $\nu^{\pm}(K)<\infty$ for any compact set $K\subset {\mathbb R}^d$ and 
there exist $\beta>0$, $c>0$, 
and $R>0$ by \eqref{eq:resolvent} and the definition of ${\cal K}_{\infty}(1)$ such that  
$$
\int_{|y|\geq R}\nu^{\pm}({\rm d}y)
=\int_{|y|\geq R}\frac{e^{-\sqrt{2}|y|}}{|y|^{(d-1)/2}}e^{\sqrt{2}|y|}|y|^{(d-1)/2}\,\nu^{\pm}({\rm d}y)
\leq c\sup_{x\in {\mathbb R}^d}\int_{|y|\geq R}G_1(x,y)\,\nu_{\beta}^{\pm}({\rm d}y)<\infty.
$$

\subsection{Population growth and spread rate on the survival event}

We first show that the results in \cite[Theorem 2.8 and Corollary 2.9]{S18} 
are valid even if $p_0\not\equiv 0$. 
For $R>0$, let $B_R:=\left\{x\in {\mathbb R}^d \mid |x|<R\right\}$ and $Z_t^{R}:=Z_t(B_R^c)$.
Define for $\delta>0$,
$$
\Lambda_{\delta}:=
\begin{cases}
\displaystyle \lambda+\sqrt{-2\lambda}\delta
& (\delta\leq \sqrt{-2\lambda}),\\
\displaystyle 
\frac{\delta^2}{2}
& (\delta>\sqrt{-2\lambda}).
\end{cases}$$

\begin{thm}\label{thm:main}
Under Assumption {\rm \ref{assum:1}}, the next assertions hold.
\begin{enumerate}
\item[{\rm (i)}] If $\delta>\sqrt{-\lambda/2}$, then 
$$\lim_{t\rightarrow\infty}Z_t^{\delta t}=0, \quad \text{${\mathbf P}_x$-a.s.}$$
\item[{\rm (ii)}] If ${\mathbf P}_x(M_{\infty}>0)>0$, 
then for any $\delta\in[0,\sqrt{-\lambda/2})$,
$$\lim_{t\rightarrow\infty}\frac{1}{t}\log Z_t^{\delta t}=-\Lambda_{\delta}, \quad \text{${\mathbf P}_x(\cdot\mid M_{\infty}>0)$-a.s.}$$
\end{enumerate}
\end{thm}

This result says that for $\delta>\sqrt{-\lambda/2}$, 
all the particles at time $t$ will be inside the ball $B_{\delta t}$ 
for all sufficiently large time $t>0$. 
On the other hand, for $\delta<\sqrt{-\lambda/2}$, 
the population outside the ball $B_{\delta t}$ at time $t$ grows exponentially with rate $-\Lambda_{\delta}$. 
\medskip

Let $L_t$ be the maximal norm of the particles alive at time $t$:
$$L_t:=
\begin{cases}
\max_{1\leq k\leq Z_t}|{\mathbf B}_t^k| & \text{($t<e_0$)},\\
0 & \text{($t\geq e_0$}).
\end{cases}$$
By the same way as for the proof of  \cite[Corollary 3.4]{S18}, 
the next corollary follows from Theorem \ref{thm:main}. 
\begin{cor}\label{cor:main}
Under Assumption {\rm \ref{assum:1}}, if ${\mathbf P}_x(M_{\infty}>0)>0$, then
$$\lim_{t\rightarrow\infty}\frac{L_t}{t}=\sqrt{-\frac{\lambda}{2}}, \quad \text{${\mathbf P}_x(\cdot \mid M_{\infty}>0)$-a.s.}$$ 
\end{cor}

\begin{rem}\label{rem:direction}\rm 
By the same way as for the proofs of Theorem \ref{thm:main} and Corollary \ref{cor:main}, 
we can show that the population growth and spread rate are uniform in direction. 
Let $\langle \cdot,\cdot\rangle$ be the standard inner product in ${\mathbb R}^d$. 
For a unit vector $r\in {\mathbb R}^d$, 
define $B_R^r:=\left\{x\in {\mathbb R}^d \mid \langle x,r\rangle<R \right\}$ 
and $Z_t^{\delta t, r}:=Z_t((B_{\delta t}^r)^c)$. Let 
$$L_t^r:=
\begin{cases}
\max_{1\leq k\leq Z_t}\langle {\mathbf B}_t^k,r\rangle & \text{($t<e_o$)},\\
0 & \text{($t\geq e_o$)}
\end{cases}$$
be the maximal displacement in direction $r$ of particles alive at time $t$. 
For $t<e_0$, we denote by $K_r(t)$ the index of a particle at time $t$ such that 
$L_t^r=\langle {\mathbf B}_t^{K_r(t)},r\rangle$. 
Then
\begin{itemize}
\item Theorem \ref{thm:main} holds 
by replacing $Z_t^{\delta t}$ with $Z_t^{\delta t,r}$;
\item Corollary \ref{cor:main} holds by replacing $L_t$ with $L_t^r$ 
and 
$$\lim_{t\rightarrow\infty}\frac{{\mathbf B}_t^{K_r(t)}}{t}
=\sqrt{-\frac{\lambda}{2}}r, \quad \text{${\mathbf P}_x(\cdot \mid M_{\infty}>0)$-a.s.}$$ 
\end{itemize}
We omit the proof of these assertions 
because the argument is similar to that of \cite[Theorem 2.10]{S18}   
by using  Remark \ref{rem:asymp-fk} and Lemma \ref{lem:lower-growth} below.
\end{rem}

\begin{rem}\label{rem:survival}\rm 
(i) \ If $\mu\in {\cal K}_{\infty}$ and  
$$\sup_{x\in {\mathbb R}^d}\sum_{n=0}^{\infty}n^2p_n(x)<\infty,$$
then $M_t\in L^2({\mathbf P}_x)$ by \cite[Lemma 3.4]{S08} 
and thus ${\mathbf P}_x(M_{\infty}>0)>0$.

\noindent
(ii) \ Suppose that $\mu\in {\cal K}_{\infty}$ and ${\mathbf P}_x(M_{\infty}>0)>0$. 
If $d=1, 2$, then by Proposition \ref{prop:equiv} below, 
we have for any $\delta>0$, 
$$\lim_{t\rightarrow\infty}\frac{1}{t}\log Z_t^{\delta t}=-\Lambda_{\delta}, \quad \text{${\mathbf P}_x(\cdot \mid e_0=\infty)$-a.s.}$$
and 
$$\lim_{t\rightarrow\infty}\frac{L_t}{t}=\sqrt{-\frac{\lambda}{2}}, \quad \text{${\mathbf P}_x(\cdot \mid e_0=\infty)$-a.s.}$$
On the other hand, if $d\geq 3$, 
then ${\mathbf P}_x(M_{\infty}=0)>0$ as in \cite[Remark 2.7]{S18}. 
Moreover, since branching occurs only finite times 
on the event $\{M_{\infty}=0\}$ by Proposition \ref{prop:equiv},  
$Z_t$ becomes a random positive constant eventually 
on the event $\{e_0=\infty\}\cap \{M_{\infty}=0\}$. 
Therefore,
$$\limsup_{t\rightarrow\infty}\frac{L_t}{\sqrt{2t\log\log t}}=1, \quad 
\text{${\mathbf P}_x(\cdot\mid \{e_0=\infty\}\cap\{M_{\infty}=0\})$-a.s.}$$    
The remarks here apply to $Z_t^{\delta t,r}$ and $L_t^r$ in Remark \ref{rem:direction}.
\end{rem}

\subsection{Upper deviation for the maximal position}   
We next show that the upper deviation of $L_t$ in \cite[Lemma 3.10]{S18} is true 
local uniformly with respect to the initial point 
and even if we allow the possibility of extinction.

\begin{thm}\label{thm:upper} 
Under Assumption {\rm \ref{assum:1}},
the next assertions hold. 
\begin{enumerate} 
\item[{\rm (i)}] If $\delta\geq \sqrt{-2\lambda}$, 
then for any compact set $K\subset {\mathbb R}^d$,
\begin{equation}\label{eq:uniform}
\lim_{t\rightarrow\infty}\frac{1}{t}\log \inf_{x\in K}{\mathbf P}_x\left(\frac{L_t}{t}\geq \delta\right)
=\lim_{t\rightarrow\infty}\frac{1}{t}\log \sup_{x\in K}{\mathbf P}_x\left(\frac{L_t}{t}\geq \delta\right)
=-\frac{\delta^2}{2}.
\end{equation}

\item[{\rm (ii)}] If $\sqrt{-\lambda/2}<\delta<\sqrt{-2\lambda}$ and ${\mathbf P}_x(M_{\infty}>0)>0$, 
then for any $x\in {\mathbb R}^d$ and for any compact set $K\subset {\mathbb R}^d$,
\begin{equation}\label{eq:uniform-1}
\lim_{t\rightarrow\infty}\frac{1}{t}\log {\mathbf P}_x\left(\frac{L_t}{t}\geq \delta\right)
=\lim_{t\rightarrow\infty}\frac{1}{t}\log \sup_{y\in K}{\mathbf P}_y\left(\frac{L_t}{t}\geq \delta\right)
=-\lambda-\sqrt{-2\lambda}\delta.
\end{equation}
\end{enumerate}
\end{thm}

Under restricted conditions,  
we can get the decay rate of ${\mathbf P}_x(L_t/t \geq \delta)$ as $t\rightarrow\infty$
more precisely. 

\begin{thm} \label{thm:decay-rate}
Assume that $p_0\equiv 0$ and 
$\mu$ is a Kato class measure with compact support in ${\mathbb R}^d$. 
If $\lambda<0$, then the next assertions hold. 
\begin{enumerate} 
\item[{\rm (i)}] If $\delta\geq \sqrt{-2\lambda}$, 
then for any $x\in {\mathbb R}^d$, 
there exist positive constants 
$C_1$, $C_2$ and $T$ such that for all $t\geq T$,
\begin{equation}\label{eq:point}
C_1e^{-\delta^2 t/2}t^{(d-2)/2} \leq 
{\mathbf P}_x\left(\frac{L_t}{t}\geq \delta\right)
\leq C_2 e^{-\delta^2 t/2}t^{(d-2)/2}.
\end{equation}

\item[{\rm (ii)}] If $\sqrt{-\lambda/2}<\delta<\sqrt{-2\lambda}$ 
and $\sup_{x\in {\mathbb R^d}}\sum_{n=1}^{\infty}n^2 p_n(x)<\infty$, 
then for any $x\in {\mathbb R}^d$, 
there exist positive constants 
$C_3$, $C_4$ and $T$ such that for all $t\geq T$, 
\begin{equation}\label{eq:point-1}
C_3e^{(-\lambda-\sqrt{-2\lambda}\delta)t}t^{(d-2)/2}\leq 
{\mathbf P}_x\left(\frac{L_t}{t}\geq \delta\right)
\leq C_4e^{(-\lambda-\sqrt{-2\lambda}\delta)t}t^{(d-1)/2}.
\end{equation}
\end{enumerate}
The lower bound of \eqref{eq:point} is valid even if $p_0\not\equiv 0$.
\end{thm}

\begin{rem}\rm 
If $\delta\geq \sqrt{-2\lambda}$, 
then by \eqref{eq:point}, \eqref{eq-br-fk} and Remark \ref{rem:asymp-fk} below,  
\begin{equation}\label{eq:equiv-p-ex}
{\mathbf P}_x\left(\frac{L_t}{t}\geq \delta\right)={\mathbf P}_x(Z_t^{\delta t}\geq 1)
\asymp{\mathbf E}_x\left[Z_t^{\delta t}\right]
\quad (t\rightarrow\infty).
\end{equation}
However, we do not know if \eqref{eq:equiv-p-ex} holds 
for $\sqrt{-\lambda/2}<\delta<\sqrt{-2\lambda}$. 

Chauvin and Rouault \cite[Theorems 2 and 3]{CR88} established 
a precise asymptotic behavior of the form like \eqref{eq:equiv-p-ex} and a Yaglom type theorem 
for branching Brownian motions on ${\mathbb R}$ with constant branching rate and mechanism. 
\end{rem}

\subsection{Population growth at the critical phase}
According to  Theorem \ref{thm:main} and \cite[Theorem 2.8]{S18}, 
the growth order of $Z_t^{\delta t}$ undergoes the phase transition at $\delta=\sqrt{-\lambda/2}$.  
We finally discuss this order at the critical phase   
under some restricted condition. 
For $\varepsilon>0$, 
let $R_{\varepsilon}(t)=t^{(d+3)/2}(\log t)(\log \log t)^{1+\varepsilon}$ 
and $r_{\varepsilon}(t)=t^{(d-2)/2}/(\log \log t)^{\varepsilon}$. 
\begin{thm}\label{thm-critical}
Assume that  $p_0\equiv 0$ and $\lambda<0$. 
Let $\delta=\sqrt{-\lambda/2}$. 
\begin{enumerate}
\item[{\rm  (i)}] 
If $\mu$ is compactly supported in ${\mathbb R}^d$, 
then for any $\varepsilon>0$, 
$${\mathbf P}_x\left(\text{there exists $T>0$ such that 
$Z_t^{\delta t}\leq R_{\varepsilon}(t)$ for all $t\geq T$}\right)=1.$$
\item[{\rm (ii)}] 
If $d\geq 3$ and ${\mathbf P}_x(M_{\infty}>0)>0$, then for any $\varepsilon>0$,
$${\mathbf P}_x\left(\text{there exists $T>0$ such that 
$Z_t^{\delta t}\geq r_{\varepsilon}(t)$ for all $t\geq T$} \mid M_{\infty}>0\right)=1.$$
\end{enumerate}
\end{thm}
If $d\geq 3$, 
then under the full conditions of Theorem \ref{thm-critical} (i) and (ii), 
$Z_t^{\delta t}$ grows polynomially ${\mathbf P}_x(\cdot\mid M_{\infty}>0)$-a.s.\
at $\delta=\sqrt{-\lambda/2}$ and the growth rate depends on the spatial dimension $d$. 
However, the exact growth rate remains unknown.

\subsection{Examples}
We apply the results in this paper to some concrete 
branching Brownian motions on ${\mathbb R}^d$. 
\begin{exam} \rm
Assume that $d=1$. 
Let $\delta_a$ be the Dirac measure at $a\in {\mathbb R}$. 
For $\gamma>0$, 
let $G_{\alpha}^{\gamma \delta_0}(x,y)$ be the $\alpha$-resolvent of the one dimensional Brownian motion 
killed by $\gamma \delta_0$:
$$G_{\alpha}^{\gamma \delta_0}(x,y)
=\frac{1}{\sqrt{2\alpha}}\left(e^{-\sqrt{2\alpha}|x-y|}-\frac{\gamma}{\sqrt{2\alpha}+\gamma}e^{-\sqrt{2\alpha}(|x|+|y|)}\right)$$
(see, e.g., \cite[p.123, 7]{BS02}).

For $a>0$, let $\mu=\gamma\delta_0-\beta\delta_a$ for some $\beta>0$ and $\gamma>0$, 
and let $\lambda=\lambda(\mu)$.
By the same way as in \cite[Example 4.1]{S08}, we have 
$$1=\beta G_{-\lambda}^{\gamma \delta_0}(a,a)
=\frac{\beta}{\sqrt{-2\lambda}}\left(1-\frac{\gamma}{\sqrt{-2\lambda}+\gamma}e^{-2a\sqrt{-2\lambda}}\right).$$
If we let $A=\sqrt{-2\lambda}$, then the equality above becomes  
$$A^2-(\beta-\gamma)A=\beta\gamma(1-e^{-2aA}).$$
This equation has a positive solution if and only if $\beta>\gamma/(1+2a\gamma)$. 
Note that this condition is derived by Takeda \cite[Example 3.10]{T02}. 

Let $\overline{\mathbf M}$ be a branching Brownian motion on ${\mathbb R}$ with branching rate 
$\mu=\delta_0+\delta_a$ and branching mechanism $\{p_n(x)\}_{n\geq 0}$ such that 
$p_0(x)+p_2(x)\equiv 1$. 
We let $p=p_2(0)$ and $q=p_2(a)$ so that  $Q(0)=2p$ and $Q(a)=2q$.  
Assume that $q\geq p$ for simplicity. 
Then $\lambda((Q-1)\mu)<0$ if one of the next conditions hold:
\begin{itemize}
\item $p\geq 1/2$ and $q>1/2$;
\item $p<1/2$, $q>1/2$ and  
$$2q-1>\frac{1-2p}{1+2a(1-2p)}.$$
\end{itemize}
In particular, Theorem \ref{thm:main}, Corollary \ref{cor:main}, Remark \ref{rem:survival} 
and Theorem \ref{thm:upper} hold under one of these conditions. 

Let $\overline{\mathbf M}$ be a branching Brownian motion on ${\mathbb R}$ with branching rate 
$\mu=c\delta_0$ for some $c>0$ and branching mechanism $\{p_n(x)\}_{n\geq 0}$ such that $p_2(0)=1$. 
Then Theorems \ref{thm:decay-rate} and \ref{thm-critical} (i) are valid 
${\mathbf P}_x$-a.s.\ with $\lambda=-c^2/2$ by Remark \ref{rem:survival}. 
Theorem \ref{thm:main} and Corollary \ref{cor:main} are proved by Bocharov and Harris \cite{BH14}.   
Corollary \ref{cor:main} also follows from \cite{E84}.
\end{exam}

\begin{exam}\rm 
Assume that $d\geq 2$. 
For $r>0$, let $\delta_r$ be the surface measure 
on the sphere $\{x\in {\mathbb R}^d \mid |x|=r\}$. 
Let $\mu=\gamma\delta_r-\beta \delta_R$  for $\beta>0$ and $\gamma>0$,  
and let $\lambda=\lambda(\mu)$. 
Define 
$$\check{\lambda}=\inf\left\{\frac{1}{2}\int_{{\mathbb R}^d}|\nabla u|^2\,{\rm d}x
+\beta\int_{{\mathbb R}^d}u^2\,{\rm d}\delta_r 
\mid u\in C_0^{\infty}({\mathbb R}^d), \gamma\int_{{\mathbb R}^d} u^2\,{\rm d}\delta_R=1\right\}.$$
Then by \cite[Lemma 2.2]{TT07}, 
$\lambda<0$ if and only if $\check{\lambda}<1$. 
\begin{enumerate}
\item[(i)] Assume first that $r<R$. 
Then by \cite[Example 3.10]{T02},  
$$\check{\lambda}=
\begin{cases}
\displaystyle \frac{\beta (r/R)}{\gamma (1+2\beta r \log(R/r))} & (d=2),\\
\displaystyle \frac{d-2}{\gamma}\left[\frac{\beta(r/R)^{d-1}}{d-2+2\beta r(1-(r/R)^{d-2})}+\frac{1}{2R}\right] & (d\geq 3).
\end{cases}$$
In particular,  $\lambda<0$ if and only if $\check{\lambda}<1$, that is, 
\begin{itemize}
\item 
$\displaystyle \gamma>\frac{\beta (r/R)}{1+2\beta r \log (R/r)}$ \ ($d=2$),
\item 
$\displaystyle \gamma>(d-2)\left[\frac{\beta(r/R)^{d-1}}{d-2+2\beta r(1-(r/R)^{d-2})}+\frac{1}{2R}\right]$ 
\ ($d\geq 3$).
\end{itemize}
\item[(ii)] Assume next that $r>R$. 
Then by the same way as in \cite[Example 3.10]{T02} together with \cite[2.3.1 in p.398--399 and 2.3.1 in p.507]{BS02}, 
we have 
$$\check{\lambda}=
\begin{cases}
\displaystyle \frac{\beta (r/R)}{\gamma (1+2\beta r \log(r/R))} & (d=2),\\
\displaystyle \frac{(d-2)(d-2+2\beta r)}{2\gamma R(d-2+2\beta r(1-(R/r)^{d-2}))} & (d\geq 3).
\end{cases}$$
Therefore, $\lambda<0$ if and only if $\check{\lambda}<1$, 
that is, 
\begin{itemize}
\item 
$\displaystyle \gamma>\frac{\beta (r/R)}{1+2\beta r \log(r/R)}$ \ ($d=2$),\\
\item 
$\displaystyle \gamma>\frac{(d-2)(d-2+2\beta r)}{2R(d-2+2\beta r(1-(R/r)^{d-2}))}$ 
\ ($d\geq 3$).
\end{itemize}
\end{enumerate}

Let $\overline{\mathbf M}$ be a branching Brownian motion on ${\mathbb R}^d$
with branching rate $\mu=\delta_r+\delta_R$ 
and spherically symmetric branching mechanism $\{p_n(x)\}_{n=0}^{\infty}$ 
such that $p_0(x)+p_2(x)\equiv 1$. 
We use the notation $p_n(x)=p_n(|x|)$.
Let $p=p_2(r)$ and $q=p_2(R)$.  
If $p<1/2$ and $q>1/2$, 
then by using (i) and (ii) with $\beta=1-2p$ and $\gamma=2q-1$, 
we can give a necessary and sufficient condition 
for $\lambda((Q-1)\mu)<0$ in terms of $p$ and $q$. 
Theorem \ref{thm:main}, Corollary \ref{cor:main}, Remark \ref{rem:survival} 
and Theorem \ref{thm:upper} are valid under this condition.

Let $\overline{\mathbf M}$ be a binary branching Brownian motion with branching rate 
$\mu=c\delta_R$ for some $c>0$. 
Then $\lambda<0$ if and only if $c>(d-2)/2$. 
Theorems \ref{thm:decay-rate} and \ref{thm-critical} are valid under this condition. 
\end{exam}

\section{Growth of Feynman-Kac functionals}\label{sect:growth}
To prove results in Section \ref{sect:result}, 
we reveal the growth rate of the expectation of $Z_t^{\delta t}$. 
By \eqref{eq-br-fk}, 
\begin{equation}\label{eq:br-fk-1}
{\mathbf E}_x\left[Z_t^{\delta t}\right]=E_x\left[e^{A_t^{\nu}};|B_t|\geq \delta t\right].
\end{equation}
Then $\nu$ is a signed measure in general because we allow $p_0\not\equiv 0$. 
In what follows, we discuss the growth rate of the expectation 
similar to that at the right hand side of \eqref{eq:br-fk-1}.

Let $\mu^+$ and $\mu^{-}$ be positive Radon measures on ${\mathbb R}^d$ in ${\cal K}_{\infty}(1)$. 
Let $\mu=\mu^+-\mu^-$ and  $\lambda:=\lambda(\mu)$. 
For $\delta>0$, we define 
$$
\Lambda_{\delta}:=
\begin{cases}
\displaystyle \lambda+\sqrt{-2\lambda}\delta
& (\delta\leq \sqrt{-2\lambda}),\\
\displaystyle 
\frac{\delta^2}{2}
& (\delta>\sqrt{-2\lambda}).
\end{cases}$$
Let $a(t)$ be a function on $(0,\infty)$ such that $a(t)=o(t) \ (t\rightarrow\infty)$ 
and $R(t):=\delta t+a(t)$.

\begin{thm}\label{thm-fk-g}
If Assumption {\rm \ref{assum:1}} is satisfied  by replacing $\nu$ with $\mu$, 
then for any $\delta>0$, $x\in {\mathbb R}^d$ and for any compact set $K\subset{\mathbb R}^d$,
$$
\lim_{t\rightarrow\infty}\frac{1}{t}\log \sup_{y\in K}E_y\left[e^{A_t^{\mu}};|B_t|\geq R(t)\right]
=\lim_{t\rightarrow\infty}\frac{1}{t}\log E_x\left[e^{A_t^{\mu}};|B_t|\geq  R(t)\right]
=-\Lambda_{\delta}.$$
\end{thm}

In \cite[Proposition 3.3]{S18}, we proved Theorem \ref{thm-fk-g} under the condition that $\mu^{-}=0$.
The proof relied on the $L^p$-independence of the spectral bounds of the Sch\"ordinger type operator (\cite{T98,T03,T08}) 
together with the fact that $A_t^{\mu}$ is nondecreasing for $\mu^{-}=0$. 
Instead of these properties, we make use of the gaugeability for Feynman-Kac semigroups developed 
by \cite{C02,T02}.

\begin{rem}\label{rem:asymp-fk}\rm
\begin{enumerate}
\item[(i)] Let $r$ be a unit vector in ${\mathbb R}^d$. 
Then the assertion in Theorem \ref{thm-fk-g} is true even if 
we replace $|B_t|$ by $\langle B_t,r\rangle$.
The proof is similar to that of Theorem \ref{thm-fk-g} 
by noting that 
$\{\langle B_t,r\rangle\}_{t\geq 0}$ is the one dimensional Brownian motion. 
\item[(ii)] Suppose that $\mu^+$ is a Kato class measure with compact support in ${\mathbb R}^d$ 
and $\mu^{-}=0$. 
Then for any $x\in {\mathbb R}^d$, we have as $t\rightarrow\infty$,
\begin{equation}\label{eq:fk-asymp-f}
E_x\left[e^{A_t^{\mu}}; |B_t|\geq R(t)\right]\asymp 
\begin{cases}
e^{(-\lambda t-\sqrt{-2\lambda}R(t))}t^{(d-1)/2} & (\delta<\sqrt{-2\lambda}), \\
e^{-R(t)^2/(2t)}t^{(d-2)/2} & (\delta\geq \sqrt{-2\lambda}).
\end{cases}
\end{equation}
In \cite[Proposition 3.1]{S18}, we proved this result  
under the condition that  $a(t)\equiv 0$, 
but the proof still works for $a(t)\ne 0$. 
If we replace $|B_t|$ by $\langle B_t,r\rangle$ in \eqref{eq:fk-asymp-f}, 
then the consequence is valid with $d=1$.
\end{enumerate}
\end{rem}

We split the proof of Theorem \ref{thm-fk-g} into the following three lemmas.

\begin{lem} \label{lem:gauge}
If $\lambda<0$, then 
there exists $p_*>1$ for any $\varepsilon>0$ such that 
for all $p\in (1,p_*)$,
$$\sup_{x\in {\mathbb R}^d}
E_x\left[\sup_{t\geq 0}\left(e^{p(\lambda-\varepsilon)t+A_t^{p\mu}}\right)\right]<\infty.$$
\end{lem}
\pf \ 
Let $\varepsilon>0$ and $p>1$. Since $\lambda-\varepsilon<0$, we have for any $t\geq 0$,
\begin{equation}\label{eq:a-f}
\begin{split}
e^{p(\lambda-\varepsilon)t+A_t^{p\mu}}
=e^{p(\lambda-\varepsilon)t}\left(1+\int_0^t e^{A_s^{p\mu}}\,{\rm d}A_s^{p\mu}\right)
&\leq e^{p(\lambda-\varepsilon)t}\left(1+\int_0^t e^{A_s^{p\mu}}\,{\rm d}A_s^{p\mu^+}\right)\\
&\leq 1+\int_0^{\infty} e^{p(\lambda-\varepsilon)u}e^{A_s^{p\mu}}\,{\rm d}A_s^{p\mu^+},
\end{split}
\end{equation}
which implies that 
\begin{equation}\label{eq:subprocess}
\begin{split}
E_x\left[\sup_{t\geq 0}\left(e^{p(\lambda-\varepsilon)t+A_t^{p\mu}}\right)\right]
&\leq 1+E_x\left[\int_0^{\infty} e^{p(\lambda-\varepsilon)s}e^{A_s^{p\mu}}\,{\rm d}A_s^{p\mu^+}\right]\\
&=1+\hat{E}_x\left[\int_0^{\zeta}e^{A_s^{p\mu}}\,{\rm d}A_s^{p\mu^+}\right].
\end{split}
\end{equation}
Here $\hat{P}_x$ and $\zeta$ are the law and life time, 
respectively, of the killed process of ${\mathbf M}$ 
by the exponential distribution with rate $p(-\lambda+\varepsilon)$.

Since 
\begin{equation*}
\begin{split}
&\inf\left\{\frac{1}{2}\int_{{\mathbb R}^d} |\nabla u|^2 \,{\rm d}x
+p(-\lambda+\varepsilon)\int_{{\mathbb R}^d}u^2\,{\rm d}x-p\int_{{\mathbb R}^d}u^2\,{\rm d}\mu 
\mid u\in C_0^{\infty}({\mathbb R}^d), \int_{{\mathbb R}^d}u^2\,{\rm d}x=1\right\}\\
&=\lambda(p\mu)+p(-\lambda+\varepsilon)\rightarrow \varepsilon>0 \ (p\rightarrow 1+0),
\end{split}
\end{equation*}
we see by \eqref{eq:subprocess}, \cite[Lemma 3.5]{T02} and \cite[Corollary 2.9 and Theorem 5.2]{C02} 
that there exists $p_*>1$ for any $\varepsilon>0$ such that for any $p\in (1,p_*)$, 
\begin{equation}\label{eq:subprocess-1}
\sup_{x\in {\mathbb R}^d}E_x\left[\sup_{t\geq 0}\left(e^{p(\lambda-\varepsilon)t+A_t^{p\mu}}\right)\right]
\leq 1+\sup_{x\in {\mathbb R}^d}\hat{E}_x\left[\int_0^{\zeta}e^{A_s^{p\mu}}\,{\rm d}A_s^{p\mu^+}\right]<\infty.
\end{equation}
This completes the proof. 
\qed

\begin{lem}\label{lem:upper-bound}
Under the same setting as in Theorem {\rm \ref{thm-fk-g}}, 
for any compact set $K\subset {\mathbb R}^d$,
\begin{equation}\label{eq:upper-bound}
\limsup_{t\rightarrow\infty}\frac{1}{t}\log \sup_{y\in K}E_y\left[e^{A_t^{\mu}};|B_t|\geq R(t)\right]
\leq -\Lambda_{\delta}.
\end{equation}
\end{lem}

\pf  \ Since 
$$e^{A_t^{\mu}}=1+\int_0^t e^{A_s^{\mu}}\,{\rm d}A_s^{\mu}\leq 1+\int_0^t e^{A_s^{\mu}}\,{\rm d}A_s^{\mu^+},$$
we have by the Markov property (see \cite[p.186, Exercise 1.13]{RY99} and \cite[(3.8)]{S18}),
\begin{equation}\label{eq:i-ii}
\begin{split}
&E_x[e^{A_t^{\mu}};|B_t|\geq R(t)]
\leq E_x\left[1+\int_0^t e^{A_s^{\mu}}\,{\rm d}A_s^{\mu^+};|B_t|\geq R(t)\right]\\
&=P_x(|B_t|\geq R(t))+E_x\left[\int_0^t e^{A_s^{\mu}}{\bf 1}_{\{|B_t|\geq R(t)\}}\,{\rm d}A_s^{\mu^+}\right]\\
&=P_x(|B_t|\geq R(t))+E_x\left[\int_0^t e^{A_s^{\mu}}P_{B_s}(|B_{t-s}|\geq R(t))\,{\rm d}A_s^{\mu^+}\right]
={\rm (I)}+{\rm (II)}.
\end{split}
\end{equation}

Let $K\subset{\mathbb R}^d$ be a compact set. 
Then there exist positive constants  $T$ and $c=c_{\delta,K}$ such that 
for all $x\in K$ and $t\geq T$,  
we have $|x|<R(t)$ and 
\begin{equation}\label{eq:prob-uniform}
\begin{split}
{\rm (I)}\leq P_x(|B_t-x|\geq R(t)-|x|)
&=\frac{\omega_d}{(2\pi)^{d/2}}\int_{(R(t)-|x|)/\sqrt{t}}^{\infty}e^{-r^2/2}r^{d-1}\,{\rm d}r\\
&\leq ce^{-R(t)^2/(2t)}\left(\frac{R(t)}{\sqrt{t}}\right)^{(d-2)/2},
\end{split}
\end{equation}
where $\omega_d$ is the area of the unit surface in ${\mathbb R}^d$.

For any $\varepsilon_1\in (0,\delta)$, we let 
\begin{equation}\label{eq:ii-1-ii-2}
\begin{split}
{\rm (II)}
&=E_x\left[\int_0^t e^{A_s^{\mu}}P_{B_s}(|B_{t-s}|\geq R(t)){\bf 1}_{\{|B_s|<\varepsilon_1 t\}}\,{\rm d}A_s^{\mu^+}\right]\\
&+E_x\left[\int_0^t e^{A_s^{\mu}}P_{B_s}(|B_{t-s}|\geq R(t)){\bf 1}_{\{|B_s|\geq \varepsilon_1 t\}}\,{\rm d}A_s^{\mu^+}\right]
={\rm (II)_1}+{\rm (II)_2}.
\end{split}
\end{equation}
Then by the same way as in \eqref{eq:prob-uniform}, 
\begin{equation}\label{eq:int-upper-0}
\begin{split}
{\rm (II)_1}
&\leq \frac{\omega_d}{(2\pi)^{d/2}}
E_x\left[\int_0^te^{A_s^{\mu}}\left(\int_{(R(t)-|B_s|)/\sqrt{t-s}}^{\infty}e^{-r^2/2}r^{d-1}\,{\rm d}r\right)
{\bf 1}_{\{|B_s|<\varepsilon_1 t\}}\,{\rm d}A_s^{\mu^+}\right]\\
&\leq \frac{\omega_d}{(2\pi)^{d/2}}
E_x\left[\int_0^te^{A_s^{\mu}}\left(\int_{(R(t)-\varepsilon_1t)/\sqrt{t-s}}^{\infty}e^{-r^2/2}r^{d-1}\,{\rm d}r\right)
\,{\rm d}A_s^{\mu^+}\right].
\end{split}
\end{equation}
For any $\varepsilon_2>0$, 
\begin{equation}\label{eq:int-upper}
\begin{split}
&\int_0^te^{A_s^{\mu}}\left(\int_{(R(t)-\varepsilon_1t)/\sqrt{t-s}}^{\infty}e^{-r^2/2}r^{d-1}\,{\rm d}r\right)
\,{\rm d}A_s^{\mu^+}\\
&\leq \sup_{0\leq s\leq t}\left(e^{(\lambda-\varepsilon_2)s+A_s^{\mu}}\right)
\int_0^te^{(-\lambda+\varepsilon_2)s}\left(\int_{(R(t)-\varepsilon_1 t)/\sqrt{t-s}}^{\infty}e^{-r^2/2}r^{d-1}\,{\rm d}r\right)
\,{\rm d}A_s^{\mu^+}.
\end{split}
\end{equation}
Since $-\lambda+\varepsilon_2>0$, we have
\begin{equation*}
\begin{split}
&\frac{\partial }{\partial s}\left[e^{(-\lambda+\varepsilon_2)s}
\left(\int_{(R(t)-\varepsilon_1 t)/\sqrt{t-s}}^{\infty}e^{-r^2/2}r^{d-1}\,{\rm d}r\right)\right]\\
&=(-\lambda+\varepsilon_2)e^{(-\lambda+\varepsilon_2)s}
\left(\int_{(R(t)-\varepsilon_1 t)/\sqrt{t-s}}^{\infty}e^{-r^2/2}r^{d-1}\,{\rm d}r\right)\\
&-e^{(-\lambda+\varepsilon_2)s}\frac{(R(t)-\varepsilon_1 t)}{2(t-s)^{3/2}}
\exp\left(-\frac{(R(t)-\varepsilon_1 t)^2}{2(t-s)}\right)
\left(\frac{R(t)-\varepsilon_1 t}{\sqrt{t-s}}\right)^{d-1}\\
&\geq 
-e^{(-\lambda+\varepsilon_2)s}\frac{(R(t)-\varepsilon_1 t)^2}{2(t-s)^{3/2}}
\exp\left(-\frac{(R(t)-\varepsilon_1 t)^2}{2(t-s)}\right)
\left(\frac{R(t)-\varepsilon_1 t}{\sqrt{t-s}}\right)^{d-1}.
\end{split}
\end{equation*}
Hence by the integration by parts formula,
\begin{equation*}
\begin{split}
&\int_0^te^{(-\lambda+\varepsilon_2)s}
\left(\int_{(R(t)-\varepsilon_1 t)/\sqrt{t-s}}^{\infty}e^{-r^2/2}r^{d-1}\,{\rm d}r\right)
\,{\rm d}A_s^{\mu^+}\\
&=-\int_0^t A_s^{\mu^+}
\frac{\partial }{\partial s}\left[e^{(-\lambda+\varepsilon_2)s}
\left(\int_{(R(t)-\varepsilon_1 t)/\sqrt{t-s}}^{\infty}e^{-r^2/2}r^{d-1}\,{\rm d}r\right)\right]\,{\rm d}s\\
&\leq 
\int_0^t A_s^{\mu^+}
e^{(-\lambda+\varepsilon_2)s}
\frac{(R(t)-\varepsilon_1 t)}{2(t-s)^{3/2}}
\exp\left(-\frac{(R(t)-\varepsilon_1 t)^2}{2(t-s)}\right)
\left(\frac{R(t)-\varepsilon_1 t}{\sqrt{t-s}}\right)^{d-1}\,{\rm d}s\\
&\leq \frac{(R(t)-\varepsilon_1 t)^d}{2} A_t^{\mu^+}
\int_0^t
\frac{e^{(-\lambda+\varepsilon_2)s}}{(t-s)^{(d+2)/2}}
\exp\left(-\frac{(R(t)-\varepsilon_1 t)^2}{2(t-s)}\right)\,{\rm d}s.
\end{split}
\end{equation*}
Combining this with \eqref{eq:int-upper-0} and \eqref{eq:int-upper}, we get 
\begin{equation}\label{eq:int-upper-1}
\begin{split}
{\rm (II)_1}
&\leq \frac{\omega_d}{2(2\pi)^{d/2}}
E_x\left[\sup_{0\leq s\leq t}\left(e^{(\lambda-\varepsilon_2)s+A_s^{\mu}}\right)A_t^{\mu^+}\right]\\
&\times (R(t)-\varepsilon_1 t)^d\int_0^t 
\frac{e^{(-\lambda+\varepsilon_2)s}}{(t-s)^{(d+2)/2}}
\exp\left(-\frac{(R(t)-\varepsilon_1 t)^2}{2(t-s)}\right)\,{\rm d}s.
\end{split}
\end{equation}

We will evaluate the integral in the right hand side of \eqref{eq:int-upper-1} in Appendix \ref{appendix:evaluate} below. 
Here we evaluate the expectation in the right hand side of \eqref{eq:int-upper-1}. 
For any $p>1$ and $q>1$ with $1/p+1/q=1$, we have by the Cauchy-Schwartz inequality, 
\begin{equation}\label{eq:c-s}
E_x\left[\sup_{0\leq s\leq t}\left(e^{(\lambda-\varepsilon_2)s+A_s^{\mu}}\right)A_t^{\mu^+}\right]
\leq E_x\left[\sup_{0\leq s\leq t}\left(e^{p(\lambda-\varepsilon_2)s+A_s^{p\mu}}\right)\right]^{1/p}
E_x\left[(A_t^{\mu^+})^q\right]^{1/q}.
\end{equation}
Then by Lemma \ref{lem:gauge}, 
there exists $p_*\in (1,\infty)$ for $\varepsilon_2>0$ 
such that for all $p\in (1,p_*)$,
$$E_x\left[\sup_{0\leq s\leq t}\left(e^{p(\lambda-\varepsilon_2)s+A_s^{p\mu}}\right)\right]
\leq \sup_{x\in {\mathbb R}^d}E_x\left[\sup_{s\geq 0}\left(e^{p(\lambda-\varepsilon_2)s+A_s^{p\mu}}\right)\right]
<\infty.$$
If we take $p\in (1,p_*)$ so that  $q>1$ is a positive integer, 
then by \cite[p.73, Corollary to Proposition 3.8]{CZ95}, 
there exist positive constants $c_1$ and $c_2$ such that for all $t\geq 0$,
$$E_x\left[(A_t^{\mu^+})^q \right]^{1/q}\leq (q!)^{1/q}\sup_{y\in {\mathbb R}^d}E_y\left[A_t^{\mu^+}\right]\leq c_1+c_2t.$$
Hence by \eqref{eq:c-s},  
\begin{equation}\label{eq:bdd}
\begin{split}
E_x\left[\sup_{0\leq s\leq t}\left(e^{(\lambda-\varepsilon_2)s+A_s^{\mu}}\right)A_t^{\mu^+}\right]
&\leq  c(\varepsilon_2)(c_1+c_2t).
\end{split}
\end{equation}
Then by \eqref{eq:int-upper-1}, \eqref{eq:int-asymp}, and \eqref{eq:bdd}, 
we have as $t\rightarrow\infty$,
\begin{equation}\label{eq:ii-1}
\begin{split}
({\rm II})_1&\lesssim 
\begin{cases}
e^{-(\delta-\varepsilon_1)^2t/2}t^{d/2} & (\delta>\sqrt{-2\lambda}),\\
e^{(-\lambda+\varepsilon_2)t}e^{-\sqrt{2(-\lambda+\varepsilon_2)}(R(t)-\varepsilon_1 t)}t^{(d+1)/2}  
& (\delta\leq \sqrt{-2\lambda}).
\end{cases}
\end{split}
\end{equation}

Define $\nu_c({\rm d}x)=e^{2c|x|/\varepsilon_1}\mu^+({\rm d}x)$ for $c>0$. 
Since $\nu_c\in {\cal K}_{\infty}(1)$ by assumption, 
it follows by \cite[Lemma 3.5]{T02} and \cite[Corollary 2.9 and Theorem 5.2]{C02} again that 
for any $c>-\lambda$,
$$\sup_{x\in {\mathbb R}^d}E_x\left[\int_0^{\infty}e^{-cs}e^{A_s^{\mu}}\,{\rm d}A_s^{\nu_c^+}\right]<\infty$$
and thus
\begin{equation}\label{eq:ii-2}
\begin{split}
({\rm II})_2
\leq E_x\left[\int_0^t e^{A_s^{\mu}}{\bf 1}_{\{|B_s|>\varepsilon_1 t\}}\,{\rm d}A_s^{\mu^+}\right]
&\leq e^{-ct}E_x\left[\int_0^{\infty}e^{-cs}e^{A_s^{\mu}}e^{2c|B_s|/\varepsilon_1}\,{\rm d}A_s^{\mu^+}\right]\\
&\leq e^{-ct}\sup_{x\in {\mathbb R}^d}E_x\left[\int_0^{\infty}e^{-cs}e^{A_s^{\mu}}\,{\rm d}A_s^{\nu_c^+}\right]
=c_3e^{-ct}.
\end{split}
\end{equation}
By taking $c>0$ large enough, 
we see by \eqref{eq:i-ii}, \eqref{eq:prob-uniform}, \eqref{eq:ii-1-ii-2}, \eqref{eq:ii-1} and \eqref{eq:ii-2} that  
\begin{equation*}
\begin{split}
&\limsup_{t\rightarrow\infty}\frac{1}{t}\log \sup_{x\in K}E_x\left[e^{A_t^{\mu}};|B_t|\geq R(t)\right]\\
&\leq 
\begin{cases}
\displaystyle -\frac{(\delta-\varepsilon_1)^2}{2}
& (\delta>\sqrt{-2\lambda}),\\
\displaystyle \left(-\lambda+\varepsilon_2-\sqrt{2(-\lambda+\varepsilon_2)}(\delta-\varepsilon_1)\right)
\vee \left(-\frac{\delta^2}{2}\right) 
& (\delta\leq \sqrt{-2\lambda}).
\end{cases}
\end{split}
\end{equation*}
Letting $\varepsilon_2\rightarrow+0$ 
and then $\varepsilon_1\rightarrow+0$, we arrive at \eqref{eq:upper-bound}.
\qed

\begin{lem}\label{lem:lower-bound}
Under the same setting as in Theorem {\rm \ref{thm-fk-g}}, 
if $\delta\geq \sqrt{-2\lambda}$, then 
for any compact set $K\subset {\mathbb R}^d$, 
\begin{equation}\label{eq:lower-bound-1}
\liminf_{t\rightarrow\infty}\frac{1}{t}\log \inf_{x\in K}E_x\left[e^{A_t^{\mu}};|B_t|\geq R(t)\right]\\
\geq -\frac{\delta^2}{2}.
\end{equation}
On the other hand, if $\delta<\sqrt{-2\lambda}$, then for any $x\in {\mathbb R}^d$,
\begin{equation}\label{eq:lower-bound}
\liminf_{t\rightarrow\infty}\frac{1}{t}\log E_x\left[e^{A_t^{\mu}};|B_t|\geq R(t)\right]\\
\geq -\lambda-\sqrt{-2\lambda}\delta.
\end{equation}

\end{lem}

\pf \ We first assume that $\delta\geq\sqrt{-2\lambda}$. 
For any $p>1$, we have by the Cauchy-Schwarz inequality,
\begin{equation}\label{eq:cs-2}
E_x\left[e^{A_t^{\mu}};|B_t|\geq R(t)\right]
\geq E_x\left[e^{-A_t^{\mu^-}};|B_t|\geq R(t)\right]
\geq \frac{P_x(|B_t|\geq R(t))^p}{E_x\left[e^{A_t^{\mu^{-}}/(p-1)};|B_t|\geq R(t)\right]^{p-1}}.
\end{equation}
Let 
$$\alpha_p=\inf\left\{\frac{1}{2}\int_{{\mathbb R}^d}|\nabla u|^2\,{\rm d}x
-\frac{1}{p-1}\int_{{\mathbb R}^d}u^2\,{\rm d}\mu^{-} 
\mid u\in C_0^{\infty}({\mathbb R}^d), \int_{{\mathbb R}^d}u^2\,{\rm d}x=1\right\}.$$
Since there exists $p_*>1$ such that 
$\sqrt{-2\lambda}\geq \sqrt{-2\alpha_{p_*}}>0$, 
it follows by Lemma \ref{lem:upper-bound} that  for any compact set $K\subset {\mathbb R}^d$,
\begin{equation}\label{eq:upper-uniform}
\limsup_{t\rightarrow\infty}\frac{1}{t}\log \sup_{x\in K}E_x\left[e^{A_t^{\mu^{-}}/(p_*-1)};|B_t|\geq R(t)\right]
\leq -\frac{\delta^2}{2}.
\end{equation}
By \cite[Appendix A]{S18}, we also have as $t\rightarrow\infty$,
\begin{equation}\label{eq:prob-est}
\begin{split}
P_x(|B_t|\geq R(t))\geq P_0(|B_t|\geq R(t))
&\sim \frac{\omega_d}{(2\pi)^{d/2}}e^{-R(t)^2/(2t)}\left(\frac{R(t)}{\sqrt{t}}\right)^{d-2}\\
&= \frac{\omega_d}{(2\pi)^{d/2}}e^{-\delta^2 t/2-\delta a(t)-a(t)^2/(2t)}
\left(\delta \sqrt{t}+\frac{a(t)}{\sqrt{t}}\right)^{d-2}. 
\end{split}
\end{equation}
By taking $p=p_*$ in \eqref{eq:cs-2}, 
we have \eqref{eq:lower-bound-1} 
by \eqref{eq:upper-uniform} and \eqref{eq:prob-est}.

We next assume that $\delta<\sqrt{-2\lambda}$. 
Fix $p\in (0,1)$ with $\delta/(1-p)\geq \sqrt{-2\lambda}$.  
Let $G$ be a relatively compact open subset in ${\mathbb R}^d$ and $K=\overline{G}$. 
Then by the Markov property, 
\begin{equation}\label{eq:domain}
\begin{split}
E_x\left[e^{A_t^{\mu}};|B_t|\geq R(t)\right]
&=E_x\left[e^{A_{pt}^{\mu}}E_{B_{pt}}\left[e^{A_{(1-p)t}^{\mu}}; |B_{(1-p)t}|\geq R(t)\right]\right]\\
&\geq E_x\left[e^{A_{pt}^{\mu}}E_{B_{pt}}\left[e^{A_{(1-p)t}^{\mu}}; |B_{(1-p)t}|\geq R(t)\right];B_{pt}\in G\right]\\
&\geq E_x\left[e^{A_{pt}^{\mu}};B_{pt}\in G\right]\inf_{y\in K}E_y\left[e^{A_{(1-p)t}^{\mu}}; |B_{(1-p)t}|\geq R(t)\right].
\end{split}
\end{equation}
Hence by \eqref{eq:lower-bound-1} and  \eqref{eq-fkf-asymp}, 
$$
\liminf_{t\rightarrow\infty}
\frac{1}{t}\log E_x\left[e^{A_t^{\mu}};|B_t|\geq R(t)\right]
\geq -\lambda p-\frac{\delta^2}{2(1-p)}.
$$
Since the right hand side above attains  the maximum value $-\lambda-\sqrt{-2\lambda}\delta$ 
at $p=1-\delta/\sqrt{-2\lambda}\in (0,1)$, 
we obtain \eqref{eq:lower-bound}. 
\qed
\medskip

Theorem \ref{thm-fk-g} is a consequence of  Lemmas \ref{lem:upper-bound} and \ref{lem:lower-bound}.

\begin{rem}\label{rem:lower-fk}\rm 
Suppose that $\mu^{-}$ is compactly supported in ${\mathbb R}^d$ 
and $\delta\geq \sqrt{-2\lambda}$. 
As in the proof of Lemma \ref{lem:lower-bound}, 
we take $p_*>1$ such that 
$\sqrt{-2\lambda}\geq \sqrt{-2\alpha_{p_*}}>0$. 
If $R(t)=\delta t+a(t)$ for $a(t)=O(1) \ (t\rightarrow\infty)$, 
then by the same calculation as in \eqref{eq:cs-2} and Remark \ref{rem:asymp-fk}, 
there exists $c>0$ such that for each $x\in{\mathbb R}^d$ and for all sufficiently large $t>0$,
$$
E_x\left[e^{A_t^{\mu}};|B_t|\geq R(t)\right]
\geq \frac{P_x(|B_t|\geq R(t))^{p_*}}{E_x\left[e^{A_t^{\mu^{-}}/(p_*-1)};|B_t|\geq R(t)\right]^{p_*-1}}
\geq ce^{-R(t)^2/(2t)}t^{(d-2)/2}.
$$
\end{rem}

\section{Proof of Theorem \ref{thm:main}}\label{sect:prf-1}
We first discuss the upper bound of $Z_t^{\delta t}$. 

\begin{lem}\label{lem:upper-growth}
Under Assumption {\rm \ref{assum:1}}, for any $\delta\geq  0$,
$$\limsup_{t\rightarrow\infty}\frac{1}{t}\log Z_t^{\delta t}\leq -\Lambda_{\delta}, \quad \text{${\mathbf P}_x$-a.s.}$$
\end{lem}

Under the additional condition that $p_0\equiv 0$, 
we proved Lemma \ref{lem:upper-growth} as \cite[Lemma 3.8]{S18} 
In the proof, we used the fact that $Z_t$ is nondecreasing,  
but this property fails for $p_0\not\equiv 0$. 
To avoid the use of this property, 
we modify the proof of \cite[Lemma 3.8]{S18} by Theorem \ref{thm-fk-g} 
together with the introduction of another branching Brownian motion   
such that it does not become extinct and its population distribution is 
comparable to the original one. 
This approach is similar to that of \cite[Subsection 4.2]{LS17+} 
(see also \cite[Section 3]{E84} for a similar discussion).

\pf \ Assume that $\delta>0$. 
For $i=1,2,\dots$ and for any $\varepsilon>0$, we have by the Chebyshev inequality,  
\begin{equation}\label{eq:chebyshev}
{\mathbf P}_x\left(\max_{i\leq s\leq i+1}Z_s^{\delta i}\geq e^{(-\Lambda_{\delta}+\varepsilon)i}\right)
\leq e^{-(-\Lambda_{\delta}+\varepsilon)i}{\mathbf E}_x\left[\max_{i\leq s\leq i+1}Z_s^{\delta i}\right].
\end{equation}
Then by the Markov property,
\begin{equation}\label{eq:max-markov}
\begin{split}
{\mathbf E}_x\left[\max_{i\leq s\leq i+1}Z_s^{\delta i}\right]
={\mathbf E}_x\left[{\mathbf E}_{{\mathbf B}_i}\left[\max_{0\leq s\leq 1}Z_s^{\delta i}\right]\right]
&\leq {\mathbf E}_x\left[\sum_{k=1}^{Z_i}{\mathbf E}_{{\mathbf B}_i^k}\left[\max_{0\leq s\leq 1}Z_s^{\delta i}\right]\right]\\
&=E_x\left[e^{A_i^{(Q-1)\mu}}{\mathbf E}_{B_i}\left[\max_{0\leq s\leq 1}Z_s^{\delta i}\right]\right].
\end{split}
\end{equation}

Let 
$\tilde{\overline{\mathbf M}}=(\{\tilde{\mathbf B}_t\}_{t\geq 0}, \{\tilde{\mathbf P}_{\bf x}\}_{{\bf x}\in {\mathbf X}})$ 
be a branching Brownian motion on ${\mathbb R}^d$ 
with branching rate $\mu$ and branching mechanism $\{q_n\}_{n=1}^{\infty}$ given by 
\begin{equation}\label{eq:modified}
q_1=p_0+p_1, \quad  q_n=p_n \ (n=2,3,\dots).
\end{equation}
Namely, for the process $\overline{\mathbf M}$, 
if a particle has no child at the splitting time, 
then we add one branching Brownian particle at the branching site.  
Hence if $\tilde{Z}_t(A)$ denotes the number of particles on a set $A\in {\mathcal B}({\mathbb R}^d)$ 
at time $t$ for $\tilde{\overline{\mathbf M}}$, 
then for any ${\mathbf x}\in {\mathbf X}$, 
\begin{equation}\label{eq:comparison}
{\mathbf P}_{\mathbf x}\left(\max_{0\leq s\leq t}Z_s(A)\geq k\right)\leq 
\tilde{\mathbf P}_{\mathbf x}\left(\max_{0\leq s\leq t}\tilde{Z}_s(A)\geq k\right)
\quad (t\geq 0, \ k=0,1,2,\dots).
\end{equation}

Let $\tilde{Z}_t=\tilde{Z}_t({\mathbb R}^d)$ and $\tilde{Z}_t^{R}=\tilde{Z}_t((B_R)^c)$. 
For $s\leq t$, let $\tilde{\mathbf B}_s^{(t),k}$ 
be the position at time $s$ of the $k$th particle alive at time $t$ 
for $\tilde{\overline{\mathbf M}}$.
Since $\tilde{Z}_t$ is nondecreasing, we have 
$$\max_{0\leq s\leq 1}\tilde{Z}_s^{\delta i}
\leq \sum_{k=1}^{\tilde{Z}_1}{\bf 1}_{\{\sup_{0\leq s\leq 1}|\tilde{\mathbf B}_s^{(1),k}|\geq \delta i\}}.$$
Then by \eqref{eq:comparison} and \cite[Lemma 3.6]{S18}, 
\begin{equation*}
\begin{split}
{\mathbf E}_x\left[\max_{0\leq s\leq 1}Z_s^{\delta i}\right]\leq
\tilde{\mathbf E}_x\left[\max_{0\leq s\leq 1}\tilde{Z}_s^{\delta i}\right]
&\leq \tilde{\mathbf E}_x\left[\sum_{k=1}^{\tilde{Z}_1}
{\bf 1}_{\{\sup_{0\leq s\leq 1}|\tilde{\mathbf B}_s^{(1),k}|\geq \delta i\}}\right]\\
&=E_x\left[e^{A_1^{(\tilde{Q}-1)\mu}};\sup_{0\leq s\leq 1}|B_s|\geq \delta i\right]
\end{split}
\end{equation*}
for 
$$\tilde{Q}(x)=\sum_{n=1}^{\infty}nq_n(x)=Q(x)+p_0(x)=1+\sum_{n=1}^{\infty}(n-1)p_n(x).$$
Therefore, for any $\alpha\in (0,\delta)$,
\begin{equation}\label{eq:est-max}
\begin{split}
&E_x\left[e^{A_i^{(Q-1)\mu}}{\mathbf E}_{B_i}\left[\max_{0\leq s\leq 1}Z_s^{\delta i}\right]\right]
\leq E_x\left[e^{A_i^{(Q-1)\mu}}E_{B_i}\left[e^{A_1^{(\tilde{Q}-1)\mu}};\sup_{0\leq s\leq 1}|B_s|\geq \delta i\right]\right]\\
&=E_x\left[e^{A_i^{(Q-1)\mu}}E_{B_i}\left[e^{A_1^{(\tilde{Q}-1)\mu}};\sup_{0\leq s\leq 1}|B_s|\geq \delta i\right];|B_i|\geq (\delta-\alpha)i\right]\\
&+E_x\left[e^{A_i^{(Q-1)\mu}}E_{B_i}\left[e^{A_1^{(\tilde{Q}-1)\mu}};\sup_{0\leq s\leq 1}|B_s|\geq \delta i\right]; 
|B_i|<(\delta-\alpha)i\right]={\rm (I)}+{\rm (II)}.
\end{split}
\end{equation}

By Theorem \ref{thm-fk-g}, 
there exists $N\geq 1$ for any $x\in {\mathbb R}^d$ and $\varepsilon>0$ 
such that for all $i\geq N$,
\begin{equation}\label{eq:(i)}
{\rm (I)}
\leq E_x\left[e^{A_i^{(Q-1)\mu}};|B_i|\geq (\delta-\alpha)i\right]
\sup_{x\in {\mathbb R}^d}E_x\left[e^{A_1^{(\tilde{Q}-1)\mu}}\right]
\leq c e^{(-\Lambda_{\delta-\alpha}+\varepsilon/2)i}.
\end{equation}
For any $x\in {\mathbb R}^d$ and $\theta>0$, since
$$P_x\left(\sup_{0\leq s\leq 1}|B_s-B_0|\geq\alpha i\right)
=P_0\left(\sup_{0\leq s\leq 1}|B_s|\geq \alpha i\right)
\leq e^{-\theta\alpha i}E_0\left[e^{\theta \sup_{0\leq s\leq 1}|B_s|}\right],$$
we have 
\begin{equation*}
\begin{split}
E_x\left[e^{A_1^{(\tilde{Q}-1)\mu}};\sup_{0\leq s\leq 1}|B_s-B_0|\geq \alpha i\right]
&\leq \sup_{x\in {\mathbb R}^d}E_x\left[e^{2A_1^{(\tilde{Q}-1)\mu}}\right]^{1/2}
P_x\left(\sup_{0\leq s\leq 1}|B_s-B_0|\geq\alpha i\right)^{1/2}\\
&\leq e^{-\theta\alpha i/2}
E_0\left[e^{\theta \sup_{0\leq s\leq 1}|B_s|}\right]^{1/2}
\sup_{x\in {\mathbb R}^d}E_x\left[e^{2A_1^{(\tilde{Q}-1)\mu}}\right]^{1/2}\\
&=c_1(\theta)e^{-\theta\alpha i/2}.
\end{split}
\end{equation*}
By \eqref{eq-fkf-asymp}, there exist $c>0$ and $N'\geq 1$ 
for any $x\in {\mathbb R}^d$ such that for any $i\geq N'$, 
\begin{equation}\label{eq:disc-asymp}
E_x\left[e^{A_i^{(Q-1)\mu}}\right]\leq ce^{-\lambda i}.
\end{equation}
Then for any $x\in {\mathbb R}^d$ and $i\geq N\vee N'$,
\begin{equation*}
\begin{split}
{\rm (II)}
\leq E_x\left[e^{A_i^{(Q-1)\mu}}E_{B_i}\left[e^{A_1^{(\tilde{Q}-1)\mu}};\sup_{0\leq s\leq 1}|B_s-B_0|\geq \alpha i\right]\right]
&\leq c_1(\theta) e^{-\theta\alpha i/2}E_x\left[e^{A_i^{(Q-1)\mu}}\right]\\
&\leq c_2(\theta)e^{(-\theta\alpha/2-\lambda)i}.
\end{split}
\end{equation*}
Combining this with \eqref{eq:(i)}, we see by \eqref{eq:chebyshev}, \eqref{eq:max-markov} and \eqref{eq:est-max} that for any $i\geq N\vee N'$, 
\begin{equation}\label{eq:up-bd}
\begin{split}
{\mathbf P}_x\left(\max_{i\leq s\leq i+1}Z_s^{\delta i}\geq e^{(-\Lambda_{\delta}+\varepsilon)i}\right)
\leq e^{-(-\Lambda_{\delta}+\varepsilon)i}{\mathbf E}_x\left[\max_{i\leq s\leq i+1}Z_s^{\delta i}\right]
&\leq e^{-(-\Lambda_{\delta}+\varepsilon)i}({\rm (I)}+{\rm (II)}) \\
&\leq c e^{-\kappa_1i}+c_2(\theta)e^{-\kappa_2 i}
\end{split}
\end{equation}
for 
$$\kappa_1=-\Lambda_{\delta}+\Lambda_{\delta-\alpha}+\frac{\varepsilon}{2}, 
\quad \kappa_2=-\Lambda_{\delta}+\varepsilon+\frac{\theta\alpha}{2}+\lambda.$$

We take $\alpha>0$ so small that $\kappa_1>0$, and then take $\theta>0$ so large that $\kappa_2>0$. 
Then, since it follows by \eqref{eq:up-bd} that 
$$\sum_{i=1}^{\infty}
{\mathbf P}_x\left(\max_{i\leq s\leq i+1}Z_s^{\delta i}\geq e^{(-\Lambda_{\delta}+\varepsilon)i}\right)
<\infty,$$
we have by the Borel-Cantelli lemma,
$${\mathbf P}_x\left(\text{$\max_{i\leq s\leq i+1}Z_s^{\delta i}\leq e^{(-\Lambda_{\delta}+\varepsilon)i}$ 
for all sufficiently large $i\geq 1$}\right)=1. $$
Therefore, for all sufficiently large $i\geq 1$ and for all $t>0$ with $i\leq t\leq i+1$,
$$Z_t^{\delta t}\leq \max_{i\leq s\leq i+1}Z_s^{\delta i}
\leq e^{(-\Lambda_{\delta}+\varepsilon)i}\leq (1\vee e^{\Lambda_{\delta}-\varepsilon})e^{(-\Lambda_{\delta}+\varepsilon)t},$$
which yields that 
$$\limsup_{t\rightarrow\infty}\frac{1}{t}\log Z_t^{\delta t}\leq -\Lambda_{\delta}+\varepsilon\rightarrow -\Lambda_{\delta} 
\quad (\varepsilon\rightarrow +0).$$
For $\delta=0$, we can show the same assertion 
by using \eqref{eq:chebyshev}, \eqref{eq:max-markov}, the first inequality in \eqref{eq:est-max}
and \eqref{eq:disc-asymp}. 
\qed
\medskip

We next discuss the lower bound of $Z_t^{\delta t}$.
\begin{lem}\label{lem:lower-growth}
Under Assumption {\rm \ref{assum:1}}, if ${\mathbf P}_x(M_{\infty}>0)>0$, 
then for any $\delta\in[0,\sqrt{-\lambda/2})$, 
$$\liminf_{t\rightarrow\infty}\frac{1}{t}\log Z_t^{\delta t}\geq -\Lambda_{\delta}, 
\quad \text{${\mathbf P}_x(\cdot \mid M_{\infty}>0)$-a.s.}$$
\end{lem}

Under the condition that $p_0\equiv 0$, 
we proved Lemma \ref{lem:lower-growth} as \cite[Lemma 3.9]{S18}.  
In the proof, we gave an asymptotic lower  bound of  the number of particles 
which are located outside the increasing balls over some time interval.
If branching occurs during this time interval, 
then we choose one of the offspring and chase its trajectory. 
However, if $p_0\not\equiv 0$, 
then particles may vanish at the splitting time. 
Here we will give an asymptotic lower  bound of the number of particles 
as we mentioned before under the condition that  no branching occurs during the time interval. 
In order to do so, we derive the locally uniform lower bound of 
the expectation in \eqref{eq:fk-n} below.
\medskip

\pf \ 
For $\delta=0$, the proof is complete by the inequality $M_t\leq \|h\|_{\infty}e^{\lambda t}Z_t$. 

In what follows, we assume that $\delta\in (0,\sqrt{-\lambda/2})$. 
Recall that ${\mathbf B}_t^k$ is the position of the $k$th particle alive at time $t$. 
At the splitting time of this particle, we choose one of its children and follow its trajectory. 
Repeating this procedure inductively,  
we can construct a trajectory starting from ${\mathbf B}_t^k$. 
We denote by ${\mathbf B}_s^{t,k}$ the position of such trajectory at time $s$ ($s\geq t$). 

For $t>s\geq 0$, let $D_{s,t}$ be the event defined by 
$$D_{s,t}:=\left\{\text{no branching occurs during the time interval $[s,t]$}\right\}.$$
Fix a constant $p\in (0,1)$ and a compact set $K\subset {\mathbb R}^d$. Then for each index $k$, 
\begin{equation*}
\begin{split}
&\left\{\text{${\mathbf B}_{np}^{np,k}\in K$, $|{\mathbf B}_s^{np,k}|>\delta s$ for all $s\in [n,n+1]$}\right\}
\cap D_{np,n+1}\\
&\supset 
\left\{\begin{minipage}[c]{90mm}
${\mathbf B}_{np}^{np,k}\in K$, 
$|{\mathbf B}_n^{np,k}|>|{\mathbf B}_n^{np,k}-{\mathbf B}_{np}^{np,k}|>\delta(n+1)+1$, \\ 
$\sup_{n\leq s\leq n+1}|{\mathbf B}_s^{np,k}-{\mathbf B}_n^{np,k}|<1$
\end{minipage}\right\}\cap D_{np,n+1}
=:E_n^k.
\end{split}
\end{equation*}
Let $x\in K$. Then for any $\varepsilon>0$ and $\alpha\in (0,\varepsilon)$, 
we  have by the Markov property,
\begin{equation}\label{eq:markov-prop}
\begin{split}
&{\mathbf P}_x\left(\sum_{k=1}^{Z_{np}}{\bf 1}_{E_n^k}\leq e^{(-\Lambda_{\delta}-\varepsilon)n}, 
Z_{np}(K)\geq e^{(-\lambda p-\alpha)n}\right)\\
&={\mathbf E}_x\left[{\mathbf P}_{{\mathbf B}_{np}}
\left(\sum_{k=1}^l {\bf 1}_{F_n^k}\leq e^{(-\Lambda_{\delta}-\varepsilon)n}\right)|_{l={Z_{np}}}; 
Z_{np}(K)\geq e^{(-\lambda p-\alpha)n}\right]
\end{split}
\end{equation}
for 
$$
F_n^k:=\left\{\text{\begin{minipage}[c]{95mm}
${\mathbf B}_0^{0,k}\in K$, 
$|{\mathbf B}_{n(1-p)}^{0,k}|>|{\mathbf B}_{n(1-p)}^{0,k}-{\mathbf B}_0^{0,k}|>\delta(n+1)+1$, \\
$\sup_{n(1-p)\leq s\leq n(1-p)+1}|{\mathbf B}_s^{0,k}-{\mathbf B}_{n(1-p)}^{0,k}|<1$\end{minipage}}\right\}
\cap D_{0,n(1-p)+1}.
$$

Let ${\mathbf x}=(x^1,\dots,x^m)$. Then by the Cauchy-Schwarz inequality, 
\begin{equation}\label{eq:indep-0}
\begin{split}
{\mathbf P}_{{\mathbf x}}
\left(\sum_{k=1}^m {\bf 1}_{F_n^k}\leq e^{(-\Lambda_{\delta}-\varepsilon)n}\right)
&={\mathbf P}_{{\mathbf x}}\left(\exp\left(-\sum_{k=1}^m{\bf 1}_{F_n^k}\right)\geq e^{-e^{(-\Lambda_{\delta}-\varepsilon)n}}\right)\\
&\leq e^{e^{(-\Lambda_{\delta}-\varepsilon)n}}{\mathbf E}_{{\mathbf x}}\left[\exp\left(-\sum_{k=1}^m{\bf 1}_{F_n^k}\right)\right].
\end{split}
\end{equation}
Since the events $F_n^k \ (1\leq k\leq m)$ are independent under ${\mathbf P}_{{\mathbf x}}$, 
we obtain by the inequality $1-x\leq e^{-x}$, 
\begin{equation}\label{eq:indep}
\begin{split}
{\mathbf E}_{{\mathbf x}}\left[\exp\left(-\sum_{k=1}^m{\bf 1}_{F_n^k}\right)\right]
=\prod_{k=1}^m{\mathbf E}_{x^k}\left[\exp\left(-{\bf 1}_{F_n^k}\right)\right]
&=
\prod_{1\leq  k\leq m, x^k\in K}\left\{1-(1-e^{-1}){\mathbf P}_{x^k}(F_n^k)\right\}\\
&\leq \prod_{1\leq k\leq m, \, x^k\in K}\exp\left(-(1-e^{-1}){\mathbf P}_{x^k}(F_n^k)\right).
\end{split}
\end{equation}

Let 
$$C_n=\left\{|B_{n(1-p)}|>|B_{n(1-p)}-B_0|>\delta(n+1)+1\right\}.$$
Then by the Markov property, we have for any $x\in K$, 
\begin{equation}\label{eq:prob-lower}
\begin{split}
{\mathbf P}_x(F_n^k)
&=E_x\left[e^{-A_{n(1-p)+1}^{\mu}};C_n \cap
\left\{\sup_{n(1-p)\leq s\leq n(1-p)+1}|B_s-B_{n(1-p)}|<1\right\}\right]\\
&=E_x\left[e^{-A_{n(1-p)}^{\mu}} 
E_{B_{n(1-p)}}\left[e^{-A_1^{\mu}};\sup_{0\leq s\leq 1}|B_s-B_0|<1\right]; 
C_n\right].
\end{split}
\end{equation}
Since it follows by  the Cauchy-Schwarz inequality that 
\begin{equation*}
\begin{split}
E_y\left[e^{-A_1^{\mu}};\sup_{0\leq s\leq 1}|B_s-B_0|<1\right]
\geq \frac{P_y(\sup_{0\leq s\leq 1}|B_s-B_0|<1)^2}{E_y\left[e^{A_1^{\mu}}\right]}
\geq  \frac{P_0(\sup_{0\leq s\leq 1}|B_s|<1)^2}{\sup_{z\in {\mathbb R}^d}E_z\left[e^{A_1^{\mu}}\right]},
\end{split}
\end{equation*}
there exists $c>0$ by \eqref{eq:prob-lower} such that 
\begin{equation}\label{eq:fk-n}
{\mathbf P}_x(F_n^k)\geq cE_x\left[e^{-A_{n(1-p)}^{\mu}}; C_n\right].
\end{equation}
By the same calculation as in \cite[p.141--142]{S18}, 
there exist $c'>0$ and $c''>0$ such that for any $x\in {\mathbb R}^d$, 
\begin{equation}\label{eq:prob-c}
P_x(C_n)\geq c'\int_{(\delta(n+1)+1)/\sqrt{n(1-p)}}^{\infty}e^{-r^2/2}r^{d-1}\,{\rm d}r
\geq c'' n^{(d-2)/2}\exp\left(-\frac{\delta^2 n}{2(1-p)}\right).
\end{equation}
Then by the Cauchy-Schwarz inequality again, we have for any $q>1$,  
\begin{equation}\label{eq:prob-lower-2}
\begin{split}
E_x\left[e^{-A_{n(1-p)}^{\mu}}; C_n\right]
&\geq \frac{P_x\left(C_n\right)^q}
{E_x\left[e^{A_{n(1-p)}^{\mu}/(q-1)};C_n\right]^{q-1}}\\
&\geq\frac{(c''n^{(d-2)/2}e^{-\delta^2 n/\{2(1-p)\}})^q}
{E_x\left[e^{A_{n(1-p)}^{\mu}/(q-1)};|B_{n(1-p)}|>\delta(n+1)+1\right]^{q-1}}.
\end{split}
\end{equation}

Let 
$$\beta_q:=\lambda\left(\frac{\mu}{q-1}\right)=\inf\left\{\frac{1}{2}\int_{{\mathbb R}^d}|\nabla u|^2\,{\rm d}x
-\frac{1}{q-1}\int_{{\mathbb R}^d}u^2\,{\rm d}\mu 
\mid u\in C_0^{\infty}({\mathbb R}^d), \int_{{\mathbb R}^d}u^2\,{\rm d}x=1\right\}.$$
Then $\sqrt{-2\lambda}\geq \sqrt{-2\beta_q}>0$ for some $q>1$. 
If we let $p=1-\delta/\sqrt{-2\lambda}$, 
then by Lemma \ref{lem:upper-bound}, there exists $N\geq 1$ for any $\varepsilon'>0$ 
such that for all $n\geq N$,
$$\sup_{y\in K}E_y\left[e^{A_{n(1-p)}^{\mu}/(q-1)};|B_{n(1-p)}|>\delta(n+1)+1\right]
\leq e^{(-\delta^2/\{2(1-p)^2\}+\varepsilon')n(1-p)}.$$
This implies that for any $x\in K$, 
the last term of \eqref{eq:prob-lower-2} is greater than 
$$\frac{(c_{\delta,p}n^{(d-2)/2}e^{-\delta^2 n/\{2(1-p)\}})^q}
{(e^{(-\delta^2/\{2(1-p)^2\}+\varepsilon')n(1-p)})^{q-1}}\\
=c_{\delta,p}^q n^{(d-2)q/2}e^{-\delta^2 n/\{2(1-p)\}}e^{-\varepsilon'(1-p)(q-1)n}
=:q_{n,p}.
$$
Hence it follows by \eqref{eq:fk-n} that 
${\mathbf P}_x(F_n^k)\geq q_{n,p}$ for all $n\geq N$. 
Since this yields that 
$$
\prod_{1\leq k\leq m, \, x^k\in K}\exp\left(-(1-e^{-1}){\mathbf P}_{x^k}(F_n^k)\right)\\
\leq \exp\left(-(1-e^{-1})q_{n,p}\sharp\{k\mid x^k\in K\}\right),
$$
we have by \eqref{eq:indep-0} and \eqref{eq:indep}, 
\begin{equation*}\label{eq:prob-lower-3}
\begin{split}
&{\mathbf E}_x\left[{\mathbf P}_{{\mathbf B}_{np}}
\left(\sum_{k=1}^l {\bf 1}_{F_n^k}\leq e^{(-\Lambda_{\delta}-\varepsilon)n}\right)|_{l={Z_{np}}}; 
Z_{np}(K)\geq e^{(-\lambda p-\alpha)n}\right]\\
&\leq {\mathbf E}_x\left[\exp\left(e^{(-\Lambda_{\delta}-\varepsilon)n}-(1-e^{-1})q_{n,p}Z_{np}(K)\right); 
Z_{np}(K)\geq e^{(-\lambda p-\alpha)n}\right]\\
&\leq \exp\left(e^{(-\Lambda_{\delta}-\varepsilon)n}-(1-e^{-1})q_{n,p}e^{(-\lambda p-\alpha)n}\right)\\
&=\exp\left(-e^{(-\Lambda_{\delta}-\alpha-\varepsilon'(1-p)(q-1))n}
\left(cn^{(d-2)q/2}-e^{(\varepsilon'(1-p)(q-1)-(\varepsilon-\alpha))n}\right)\right).
\end{split}
\end{equation*}
If we take $\varepsilon'>0$ so small that 
$\varepsilon'(1-p)(q-1)<\varepsilon-\alpha$, 
then there exists $N^{''}\geq 1$ such that for all $n\geq N^{''}$, 
we obtain by \eqref{eq:markov-prop}, 
\begin{equation*}
\begin{split}
&{\mathbf P}_x\left(\sum_{k=1}^{Z_{np}}{\bf 1}_{E_n^k}\leq e^{(-\Lambda_{\delta}-\varepsilon)n}, 
Z_{np}(K)\geq e^{(-\lambda p-\alpha)n}\right)
\leq \exp\left(-e^{(-\Lambda_{\delta}-\alpha-\varepsilon'(1-p)(q-1))n}
c'n^{(d-2)q/2}\right). 
\end{split}
\end{equation*}
Noting that  $\Lambda_{\delta}<0$ for any $\delta\in (0,\sqrt{-\lambda/2})$, 
we can take $\varepsilon>0$ such that $-\Lambda_{\delta}-\varepsilon>0$. 
Since this implies that 
$$\sum_{n=1}^{\infty}{\mathbf P}_x\left(\sum_{k=1}^{Z_{np}}{\bf 1}_{E_n^k}\leq e^{(-\Lambda_{\delta}-\varepsilon)n}, 
Z_{np}(K)\geq e^{(-\lambda p-\alpha)n}\right)<\infty,$$
we see by the Borel-Cantelli lemma that the event 
$$\left\{\sum_{k=1}^{Z_{np}}{\bf 1}_{E_n^k}> e^{(-\Lambda_{\delta}-\varepsilon)n}\right\}\cup 
\left\{Z_{np}(K)< e^{(-\lambda p-\alpha)n}\right\}$$
occurs infinitely often. 
By \eqref{eq-limit-thm}, we further obtain 
$${\mathbf P}_x\left(\text{$\sum_{k=1}^{Z_{np}}{\bf 1}_{E_n^k}>e^{(-\Lambda_{\delta}-\varepsilon)n}$ 
for all sufficiently large $n$} \mid M_{\infty}>0\right)=1.$$
Hence we have ${\mathbf P}_x(\cdot\mid M_{\infty}>0)$-a.s.\ for all sufficiently large $t>0$, 
\begin{equation*}
\begin{split}
Z_t^{\delta t}
=\sum_{k=1}^{Z_t}{\bf 1}_{\{|{\mathbf B}_t^k|>\delta t\}}
&\geq \sum_{k=1}^{Z_{[t]p}}
{\bf 1}_
{
\left\{\left\{
\text{${\mathbf B}_{[t]p}^{[t]p,k}\in K$, 
$|{\mathbf B}_s^{[t]p,k}|>\delta s$ for all $s\in [[t],[t]+1]$}
\right\}\cap D_{[t]p,[t]+1}\right\}}\\
&\geq \sum_{k=1}^{Z_{[t]p}}{\bf 1}_{E_{[t]}^k}\geq e^{(-\Lambda_{\delta}-\varepsilon)[t]},
\end{split}
\end{equation*}
which implies that 
$$\liminf_{t\rightarrow\infty}\frac{1}{t}\log Z_t^{\delta t}\geq -\Lambda_{\delta}-\varepsilon.$$
By letting $\varepsilon'\rightarrow +0$, $\alpha\rightarrow+0$, and then $\varepsilon\rightarrow +0$, 
we arrive at the desired conclusion.
\qed
\medskip

Theorem \ref{thm:main} is a consequence of Lemmas \ref{lem:upper-growth} and \ref{lem:lower-growth}.

\section{Proofs of Theorems \ref{thm:upper} and \ref{thm:decay-rate}}
In this section, we prove Theorems \ref{thm:upper} and \ref{thm:decay-rate}. 
For the lower bound of \eqref{eq:uniform-1} especially, 
we take into consideration the effect of $p_0\not\equiv 0$ 
as in  Lemma \ref{lem:lower-growth}.
\subsection{Proof of Theorem \ref{thm:upper}} 
By the Chebyshev inequality and \eqref{eq-br-fk},
$${\mathbf P}_x\left(\frac{L_t}{t}\geq \delta\right)
={\mathbf P}_x(Z_t^{\delta t}\geq 1)\leq {\mathbf E}_x\left[Z_t^{\delta t}\right]
=E_x\left[e^{A_t^{(Q-1)\mu}};|B_t|\geq \delta t\right]$$
for any $x\in {\mathbb R}^d$. 
Theorem \ref{thm-fk-g} then implies that for any compact set $K\subset {\mathbb R}^d$,
\begin{equation}\label{eq:upper-devi-1}
\limsup_{t\rightarrow\infty}\frac{1}{t}\log \sup_{x\in K}{\mathbf P}_x\left(\frac{L_t}{t}\geq \delta\right)
\leq \limsup_{t\rightarrow\infty}\frac{1}{t}\log \sup_{x\in K}E_x\left[e^{A_t^{(Q-1)\mu}};|B_t|\geq \delta t\right]
=-\Lambda_{\delta}.
\end{equation}

We now assume that $\delta\geq \sqrt{-2\lambda}$. Since  
$${\mathbf P}_x\left(\frac{L_t}{t}\geq \delta\right)
\geq {\mathbf P}_x\left(\frac{L_t}{t}\geq \delta, t<T\right)
=E_x\left[e^{-A_t^{\mu}}; |B_t|\geq \delta t\right],$$
Lemma \ref{lem:lower-bound} yields that for any compact set $K\subset {\mathbb R}^d$,
$$
\liminf_{t\rightarrow\infty}\frac{1}{t}\log \inf_{x\in K}{\mathbf P}_x\left(\frac{L_t}{t}\geq \delta\right)
\geq \liminf_{t\rightarrow\infty}\frac{1}{t}\log \inf_{x\in K}E_x\left[e^{-A_t^{\mu}}; |B_t|\geq \delta t\right]=-\frac{\delta^2}{2}.
$$
Combining this with \eqref{eq:upper-devi-1}, we complete  the proof of (i).

We next assume that $\sqrt{-\lambda/2}<\delta<\sqrt{-2\lambda}$. 
Fix $p\in (0,1)$ and $\alpha>0$. For $x\in {\mathbb R}^d$, 
we take a compact set $K\subset {\mathbb R}^d$ so  that $x\in K$. 
Let 
$$C_t:=\left\{Z_t(K)\geq e^{(-\lambda-\alpha)t}\right\}$$
and 
$$D_{s,t}:=\left\{\text{no branching occurs during the time interval $[s,t]$}\right\}.$$
Then by the Markov property,
\begin{equation}\label{eq:lower-deviation}
\begin{split}
{\mathbf P}_x(Z_t^{\delta t}\geq 1)
&={\mathbf P}_x\left(\bigcup_{k=1}^{Z_t}\left\{|{\mathbf B}_t^k|\geq \delta t\right\}\right)
\geq {\mathbf P}_x\left(\bigcup_{k=1}^{Z_{pt}}
\left\{{\mathbf B}_{pt}^{pt,k}\in K,  
|{\mathbf B}_t^{pt,k}|\geq \delta t\right\}, 
C_{pt}\cap D_{pt,t}\right)\\
&={\mathbf E}_x\left[{\mathbf P}_{{\mathbf B}_{pt}}
\left(\bigcup_{k=1}^{l}
\left\{{\mathbf B}_0^{0,k}\in K, |{\mathbf B}_{(1-p)t}^{0,k}|\geq \delta t\right\}\cap D_{0,(1-p)t}
\right)\mid_{l=Z_{pt}}; C_{pt}\right].
\end{split}
\end{equation}
For ${\mathbf x}=(x^1,\dots, x^l)\in {\mathbf X}$, since
\begin{equation*}
\begin{split}
&{\mathbf P}_{{\mathbf x}}
\left(\bigcup_{k=1}^{l}\left\{{\mathbf B}_0^{0,k}\in K, |{\mathbf B}_{(1-p)t}^{0,k}|\geq \delta t\right\}
\cap D_{0,(1-p)t}\right)\\
&=1-{\mathbf P}_{{\mathbf x}}\left(\bigcap_{k=1}^{l}
\left\{\left\{{\mathbf B}_0^{0,k}\in K, |{\mathbf B}_{(1-p)t}^{0,k}|\geq \delta t\right\}\cap D_{0,(1-p)t}
\right\}^c\right)\\
&=1-\prod_{k=1}^l{\mathbf P}_{x^k}
\left(
\left\{{\mathbf B}_0^{1}\in K, |{\mathbf B}_{(1-p)t}^{1}|\geq \delta t, 
T>(1-p)t\right\}^c
\right)
\\
&=1-\prod_{k=1}^l
\left(1-{\mathbf P}_{x^k}
\left({\mathbf B}_0^{1}\in K, |{\mathbf B}_{(1-p)t}^{1}|\geq \delta t,  
T>(1-p)t
\right)\right),
\end{split}
\end{equation*}
we have by \eqref{eq:lower-deviation}, 
\begin{equation}\label{eq:lower-deviation-0}
{\mathbf P}_x(Z_t^{\delta t}\geq 1)
\geq {\mathbf E}_x\left[1-
\prod_{k=1}^{Z_{pt}}
\left\{1-{\mathbf P}_{{\mathbf B}_{pt}^k}
\left({\mathbf B}_0^{1}\in K, |{\mathbf B}_{(1-p)t}^{1}|\geq \delta t, 
T>(1-p)t
\right)\right\}; C_{pt}\right].
\end{equation}

In what follows, we fix $p\in (0,1)$ with $\delta/(1-p)>\sqrt{-2\lambda}$.
Then by Lemma \ref{lem:lower-bound}, 
there exists $t_0>0$ for any $\varepsilon>0$ such that for any $y\in K$ and $t\geq t_0$,
\begin{equation*}
\begin{split}
&{\mathbf P}_y\left(|{\mathbf B}_{(1-p)t}^{1}|\geq \delta t, T>(1-p)t \right)
=E_y\left[e^{-A_{(1-p)t}^{\mu}};|B_{(1-p)t}|\geq \delta t \right]\\
&\geq \inf_{z\in K}E_z\left[e^{-A_{(1-p)t}^{\mu}};|B_{(1-p)t}|\geq \frac{\delta}{1-p} \cdot (1-p)t \right]
\geq \exp\left(\left(-\frac{\delta^2}{2(1-p)^2}-\varepsilon \right)(1-p)t\right).
\end{split}
\end{equation*}
Hence for any $t\geq t_0$,
\begin{equation}\label{eq:lower-cal}
\begin{split}
{\mathbf P}_x(Z_t^{\delta t}\geq 1)
&\geq
{\mathbf E}_x\left[1-\prod_{1\leq k\leq Z_{p t}, {\mathbf B}_{p t}^k\in K}
\left\{1-\exp\left(\left(-\frac{\delta^2}{2(1-p)^2}-\varepsilon \right)(1-p)t\right)\right\}; 
C_{pt}\right]\\
&={\mathbf E}_x\left[1-
\left\{1-\exp\left(\left(-\frac{\delta^2}{2(1-p)^2}-\varepsilon \right)(1-p)t\right)\right\}^{Z_{p t}(K)}; 
C_{pt}\right]\\
&\geq 
\left\{1-
\left(1-\exp\left(\left(-\frac{\delta^2}{2(1-p)^2}-\varepsilon \right)(1-p)t\right)\right)^{e^{(-\lambda -\alpha)pt}}\right\}
{\mathbf P}_x\left(C_{pt}\right).
\end{split}
\end{equation}

By elementary calculation as in the proof of \cite[Lemma 3.10]{S08},
\begin{equation*}
\begin{split}
&\liminf_{t\rightarrow\infty}\frac{1}{t}\log \left\{1-
\left(1-\exp\left(\left(-\frac{\delta^2}{2(1-p)^2}-\varepsilon \right)(1-p)t\right)\right)^{e^{(-\lambda-\alpha)pt}}\right\}\\
&\geq -\lambda p-\frac{\delta^2}{2(1-p)}-\varepsilon(1-p)-\alpha p \rightarrow -\lambda p-\frac{\delta^2}{2(1-p)}
\quad (\varepsilon\rightarrow +0, \ \alpha\rightarrow +0).
\end{split}
\end{equation*}
The right hand side above takes the maximal value  $-\lambda-\sqrt{-2\lambda}\delta$ at $p=1-\delta/\sqrt{-2\lambda}$.
Since it follows by \eqref{eq-limit-thm} that 
$${\mathbf P}_x(C_{pt})
\geq {\mathbf P}_x(C_{pt}\cap \{M_{\infty}>0\})\rightarrow {\mathbf P}_x(M_{\infty}>0) 
\quad (t\rightarrow\infty),$$
we see by \eqref{eq:lower-cal} that  if ${\mathbf P}_x(M_{\infty}>0)>0$ and $\delta>\sqrt{-\lambda/2}$, then 
$$\liminf_{t\rightarrow\infty}\frac{1}{t}\log {\mathbf P}_x\left(\frac{L_t}{t}\geq \delta\right)
=\liminf_{t\rightarrow\infty}\frac{1}{t}\log {\mathbf P}_x(Z_t^{\delta t}\geq 1)
\geq -\lambda-\sqrt{-2\lambda}\delta.
$$
Combining this with \eqref{eq:upper-devi-1}, we finish the proof.  
\qed

\subsection{Proof of Theorem \ref{thm:decay-rate}}

We first show (i). Assume that $\delta\geq \sqrt{-2\lambda}$. 
Since it follows by the proof of Theorem \ref{thm:upper} (i) that
\begin{equation}\label{eq:cheb}
E_x\left[e^{-A_t^{\mu}};|B_t|\geq \delta t\right]\leq 
{\mathbf P}_x\left(\frac{L_t}{t}\geq \delta\right)
\leq E_x\left[e^{A_t^{(Q-1)\mu}};|B_t|\geq \delta t\right],
\end{equation}
we get \eqref{eq:point} by Remarks \ref{rem:asymp-fk} and \ref{rem:lower-fk}.

We next show (ii). 
Assume that  $\delta\in (\sqrt{-\lambda/2},\sqrt{-2\lambda})$. 
If $p_0\equiv 0$ and $\mu$ is compactly supported in ${\mathbb R}^d$,  
then the upper bound of \eqref{eq:point-1} follows by \eqref{eq:cheb} and Remark \ref{rem:asymp-fk}. 
For the lower bound of it, we make use of the Feynman-Kac expression of ${\mathbf P}_x(L_t/t\geq \delta)$.  
Such an approach is similar to that of \cite{CR88} and due to McKean \cite{Mc75, Mc76} 
(see also \cite[Section 1.3]{INW68-1} and \cite[Example 3.4]{INW68-2}). 

Let us derive the Feynman-Kac expression of ${\mathbf P}_x(L_t/t\geq \delta)$. 
Let $f$ be a nonnegative Borel measurable function on ${\mathbb R}^d$ 
such that $0\leq f(x)\leq 1$ for any $x\in {\mathbb R}^d$ and  let
$$u(t,x)={\mathbf E}_x\left[\prod_{k=1}^{Z_t}f({\mathbf B}_t^k)\right]$$
for $t\geq 0$ and $x\in {\mathbb R}^d$. 
We first give the Feynman-Kac expression of $u(t,x)$ and $1-u(t,x)$. 
Let 
$$F_u(t,x)=\sum_{n=1}^{\infty}p_n(x)u(t,x)^n, \quad G_u(t,x)=\sum_{n=1}^{\infty}p_n(x)u(t,x)^{n-1}.$$
If we define 
$$H_u(t,x)=\sum_{n=1}^{\infty}p_n(x)\left(\sum_{k=1}^n u(t,x)^{k-1}\right),$$
then
\begin{equation}\label{eq:f-h}
(1-u(t,x))H_u(t,x)=1-p_0-F_u(t,x).
\end{equation}

\begin{lem} 
Assume that $0\leq f(x)\leq 1$ for any $x\in {\mathbb R}^d$. 
Then 
\begin{equation}\label{eq:fk-0}
u(t,x)=E_x\left[\exp\left(\int_0^t (G_u(t-s,B_s)-1)\,{\rm d}A_s^{\mu}\right)f(B_t)\right]
\end{equation}
and 
\begin{equation}\label{eq:fk-1}
\begin{split}
1-u(t,x)
&=E_x\left[\exp\left(\int_0^t (H_u(t-s,B_s)-1)\,{\rm d}A_s^{\mu}\right)(1-f(B_t))\right]\\
&+E_x\left[\int_0^t \exp\left(\int_0^s (H_u(t-w,B_w)-1)\,{\rm d}A_w^{\mu}\right)p_0(B_s)\,{\rm d}A_s^{\mu}\right].
\end{split}
\end{equation}
\end{lem}

\pf \ Let $F=F_u$, $G=G_u$ and $H=H_u$. 
For $s,t\geq 0$ with $t\geq s$, define 
$$C_s^{G,t}=\int_0^s G(t-r,B_r)\,{\rm d}A_r^{\mu}.$$
We first prove by induction that for any $n\geq 0$,
\begin{equation}\label{eq:n-n}
u(t,x)
=\sum_{k=0}^n \frac{1}{k!}E_x\left[e^{-A_t^{\mu}}f(B_t) (C_t^{G,t})^k\right]
+\frac{1}{n!}E_x\left[\int_0^t e^{-A_s^{\mu}}F(t-s,B_s)(C_s^{G,t})^n \,{\rm d}A_s^{\mu}\right].
\end{equation}
For $n=0$, this equality is valid because we have by the strong Markov property,
\begin{equation}\label{eq:n-0}
\begin{split}
u(t,x)
&={\mathbf E}_x\left[f({\mathbf B}_t^1):t<T\right]+{\mathbf E}_x\left[\prod_{k=1}^{Z_t}f({\mathbf B}_t^k);T\leq t\right]\\
&={\mathbf E}_x\left[f({\mathbf B}_t^1):t<T\right]
+{\mathbf E}_x\left[{\mathbf E}_{{\mathbf B}_T}\left[\prod_{k=1}^{Z_{t-s}}f({\mathbf B}_{t-s}^k)\right]|_{s=T};T\leq t\right]\\
&=E_x\left[e^{-A_t^{\mu}}f(B_t)\right]
+E_x\left[\int_0^t e^{-A_s^{\mu}}\sum_{n=1}^{\infty}p_n(B_s)u(t-s,B_s)^n\,{\rm d}A_s^{\mu}\right]\\
&=E_x\left[e^{-A_t^{\mu}}f(B_t)\right]+E_x\left[\int_0^t e^{-A_s^{\mu}}F(t-s,B_s)\,{\rm d}A_s^{\mu}\right].
\end{split}
\end{equation}
Suppose that \eqref{eq:n-n} is true for some $n\geq 1$. Then by \eqref{eq:n-0},
\begin{equation*}
\begin{split}
&E_x\left[\int_0^t e^{-A_s^{\mu}}F(t-s,B_s)\left(C_s^{G,t}\right)^n\,{\rm d}A_s^{\mu}\right]\\
&=E_x\left[\int_0^t e^{-A_s^{\mu}}u(t-s,B_s)G(t-s,B_s)\left(C_s^{G,t}\right)^n\,{\rm d}A_s^{\mu}\right]\\
&=E_x\left[\int_0^t e^{-A_s^{\mu}}
E_{B_s}\left[e^{-A_{t-s}^{\mu}}f(B_{t-s})\right]\left(C_s^{G,t}\right)^n\,{\rm d}C_s^{G,t}\right]\\
&+E_x\left[\int_0^te^{-A_s^{\mu}}
E_{B_s}\left[\int_0^{t-s} e^{-A_w^{\mu}}F(t-s-w,B_w)\,{\rm d}A_w^{\mu}\right]
\left(C_s^{G,t}\right)^n\,{\rm d}C_s^{G,t}\right]\\
&={\rm (I)}+{\rm (II)}.
\end{split}
\end{equation*}
By the Markov property,
\begin{equation*}
\begin{split}
{\rm (I)}&=E_x\left[\int_0^t e^{-A_s^{\mu}}
E_x\left[ e^{-A_{t-s}^{\mu}\circ\theta_s} f(B_{t-s}\circ\theta_s)\mid{\cal F}_s\right]
\left(C_s^{G,t}\right)^n\,{\rm d}C_s^{G,t}\right]\\
&=E_x\left[\int_0^t e^{-A_s^{\mu}}e^{-A_{t-s}^{\mu}\circ\theta_s} f(B_{t-s}\circ\theta_s)
\left(C_s^{G,t}\right)^n\,{\rm d}C_s^{G,t}\right]\\
&=E_x\left[e^{-A_t^{\mu}}f(B_t)\int_0^t \left(C_s^{G,t}\right)^n\,{\rm d}C_s^{G,t}\right]
=\frac{1}{n+1}E_x\left[e^{-A_t^{\mu}}f(B_t)
\left(C_t^{G,t}\right)^{n+1}\right]
\end{split}
\end{equation*}
and
\begin{equation*}
\begin{split}
{\rm (II)}
&=E_x\left[\int_0^t e^{-A_s^{\mu}}
E_x\left[\int_0^{t-s} e^{-A_w^{\mu}\circ \theta_s}F(t-s-w,B_w\circ\theta_s)
\,{\rm d}A_w^{\mu}\circ\theta_s \mid{\cal F}_s\right]
\left(C_s^{G,t}\right)^n\,{\rm d}C_s^{G,t}\right]\\
&=E_x\left[\int_0^t 
\left(\int_s^t e^{-A_w^{\mu}}F(t-w,B_w)
\,{\rm d}A_w^{\mu}\right)
\left(C_s^{G,t}\right)^n\,{\rm d}C_s^{G,t}\right]\\
&=E_x\left[\int_0^t e^{-A_w^{\mu}}F(t-w,B_w)
\left(\int_0^w \left(C_s^{G,t}\right)^n\,{\rm d}C_s^{G,t}\right)\,{\rm d}A_w^{\mu}
\right]\\
&=\frac{1}{n+1}E_x\left[\int_0^t e^{-A_w^{\mu}}F(t-w,B_w)
 \left(C_w^{G,t}\right)^{n+1}\,{\rm d}A_w^{\mu}\right].
\end{split}
\end{equation*}
Hence the induction is complete by \eqref{eq:n-n}. 

We next show that 
\begin{equation}\label{eq:n-lim}
\lim_{n\rightarrow\infty}\frac{1}{n!}E_x\left[\int_0^t e^{-A_s^{\mu}}
F(t-s,B_s)\left(C_s^{G,t}\right)^n\,{\rm d}A_s^{\mu}\right]=0.
\end{equation}
Since $0\leq u(t,x)\leq 1$, 
we have 
$$F(t,x)\leq G(t,x)\leq 1 \quad \text{for all $t>0$ and $x\in {\mathbb R}^d$,}$$
and therefore
$$\int_0^t e^{-A_s^{\mu}}F(t-s,B_s)\left(C_s^{G,t}\right)^n\,{\rm d}A_s^{\mu}
\leq \int_0^t e^{-A_s^{\mu}}(A_s^{\mu})^n \,{\rm d}A_s^{\mu}.$$
Since this implies that 
\begin{equation*}
\begin{split}
\sum_{n=0}^{\infty}\frac{1}{n!}E_x\left[\int_0^t e^{-A_s^{\mu}}F(t-s,B_s)
\left(C_s^{G,t}\right)^n\,{\rm d}A_s^{\mu}\right]
&\leq \sum_{n=0}^{\infty}\frac{1}{n!}E_x\left[\int_0^t e^{-A_s^{\mu}}(A_s^{\mu})^n \,{\rm d}A_s^{\mu}\right]\\
&=E_x\left[A_t^{\mu}\right]<\infty,
\end{split}
\end{equation*}
we get \eqref{eq:n-lim}. 
Furthermore, we obtain \eqref{eq:fk-0}  
by letting  $n\rightarrow\infty$ in \eqref{eq:n-n}.

We let $v(t,x)=1-u(t,x)$. 
Since 
$$e^{-A_t^{\mu}}=1-\int_0^te^{-A_s^{\mu}}\,{\rm d}A_s^{\mu},$$
we have by \eqref{eq:f-h} and \eqref{eq:n-0},
\begin{equation}
\begin{split}
v(t,x)
&=E_x\left[e^{-A_t^{\mu}}(1-f(B_t))\right]
+E_x\left[\int_0^te^{-A_s^{\mu}}(1-F(t-s,B_s))\,{\rm d}A_s^{\mu}\right]\\
&=E_x\left[e^{-A_t^{\mu}}(1-f(B_t))\right]
+E_x\left[\int_0^te^{-A_s^{\mu}}v(t-s,B_s)H(t-s,B_s)\,{\rm d}A_s^{\mu}\right]\\
&+E_x\left[\int_0^t e^{-A_s^{\mu}}p_0(B_s)\,{\rm d}A_s^{\mu}\right].
\end{split}
\end{equation}
Then the proof is complete by the induction and calculation similar 
to those for \eqref{eq:fk-0}.
\qed
\medskip

Let $f_R(x)={\bf 1}_{\{|x|<R\}}$ for $R>0$. 
If we define 
$$u_R(t,x)={\mathbf E}_x\left[\prod_{k=1}^{Z_t} f_R({\mathbf B}_t^k)\right]$$
and $v_R(t,x)=1-u_R(t,x)$, then $v_R(t,x)={\mathbf P}_x(L_t\geq R)$. 
We also define  
$$C_s^{R,t}=\int_0^s (H_{u_R}(t-w,B_w)-1)\,{\rm d}A_w^{\mu}$$ 
for $s,t\geq 0$ with $t\geq s$. Then 
$$C_t^{R,t}=C_s^{R,t}+C_{t-s}^{R,t}\circ\theta_s.$$
For $\delta>0$, we let $D_s^{t}=C_s^{\delta t,t}$.
Since 
$$v_{\delta t}(t,x)=1-u_{\delta t}(t,x)={\mathbf P}_x(L_t/t\geq \delta),$$ 
we have by \eqref{eq:fk-1},
\begin{equation}\label{eq:fk-ex}
{\mathbf P}_x(L_t/t\geq \delta)
=E_x\left[e^{D_t^t};|B_t|\geq \delta t\right]
+E_x\left[\int_0^t e^{D_s^t}p_0(B_s)\,{\rm d}A_s^{\mu}\right]\geq E_x\left[e^{D_t^t};|B_t|\geq \delta t\right].
\end{equation}
To derive the decay rate  of the right hand side above as $t\rightarrow\infty$, we show

\begin{lem}\label{lem:fk-d}
Suppose that $\mu$ is compactly supported in ${\mathbb R}^d$ 
and $\sup_{x\in {\mathbb R}^d}\sum_{n=1}^{\infty}n^2p_n(x)<\infty$. 
Then for any $p\in (0,1)$ and $\delta>\sqrt{-\lambda/2}$,
$$
\lim_{t\rightarrow\infty}e^{\lambda pt}E_x\left[e^{D_{pt}^{t}}\right]
=h(x)\int_{{\mathbb R}^d}h(y)\,{\rm d}y.
$$
\end{lem}

\pf \ 
For any $v>0$,  
$$
n-\sum_{k=1}^n(1-v)^{k-1}
=n-\frac{1-(1-v)^n}{v}\\
=\frac{(1-v)^n-(1-nv)}{v}\leq \frac{n(n-1)}{2}v.
$$
This inequality is true also for $v=0$. 
Therefore,
\begin{equation*}
\begin{split}
Q(x)-H_{1-v_R}(t,x)
&=\sum_{n=1}^{\infty}p_n(x)\left(n-\sum_{k=1}^n(1-v_R(t,x))^{k-1}\right)\\
&\leq \frac{1}{2}\sum_{n=1}^{\infty}n(n-1)p_n(x)v_R(t,x).
\end{split}
\end{equation*}
Then by the inequality $1-e^{-x}\leq x$, we have for any $p\in (0,1)$ and $t\geq 0$,
\begin{equation*}
\begin{split}
&1-e^{-(A_{pt}^{(Q-1)\mu}-C_{pt}^{R,t})}
=1-\exp\left(-\int_0^{pt} (Q(B_s)-H_{1-v_R}(t-s,B_s))\,{\rm d}A_s^{\mu}\right)\\
&\leq \int_0^{pt} (Q(B_s)-H_{1-v_R}(t-s,B_s))\,{\rm d}A_s^{\mu}\\
&\leq \int_0^{pt} \frac{1}{2}\sum_{n=1}^{\infty}n(n-1)p_n(B_s)v_R(t-s,B_s)\,{\rm d}A_s^{\mu}
=\int_0^{pt} v_R(t-s,B_s)\,{\rm d}A_s^{M\mu}
\end{split}
\end{equation*}
for 
$$M(x)=\frac{1}{2}\sum_{n=1}^{\infty}n(n-1)p_n(x).$$
Hence 
\begin{equation}\label{eq:dist}
0\leq e^{A_{pt}^{(Q-1)\mu}}-e^{C_{pt}^{R,t}}
=e^{A_{pt}^{(Q-1)\mu}}\left(1-e^{-(A_{pt}^{(Q-1)\mu}-C_{pt}^{R,t})}\right)
\leq e^{A_{pt}^{(Q-1)\mu}}\int_0^{pt} v_R(t-s,B_s)\,{\rm d}A_s^{M\mu}.
\end{equation}

If we take $R=\delta t$, then for any $s\in [0,pt]$,
$$
v_{\delta t}(t-s,x)
={\mathbf P}_x(L_{t-s}\geq \delta t)
\leq {\mathbf P}_x(L_{t-s}\geq \delta(t-s))
={\mathbf P}_x(Z_{t-s}^{\delta(t-s)}\geq 1)
\leq {\mathbf E}_x\left[Z_{t-s}^{\delta(t-s)}\right].
$$
Since $\delta>\sqrt{-\lambda/2}$, 
Theorem \ref{thm-fk-g} yields that 
for any compact set $K\subset {\mathbb R}^d$ and for any $\varepsilon\in (0,\Lambda_{\delta})$, there exists $T>0$ such that 
for all $t\geq T$ and $s\in[0,pt]$,
$$\sup_{x\in K}v_{\delta t}(t-s,x)\leq \sup_{x\in K}{\mathbf E}_x\left[Z_{t-s}^{\delta(t-s)}\right]
\leq e^{(-\Lambda_{\delta}+\varepsilon)(t-s)}
\leq e^{(-\Lambda_{\delta}+\varepsilon)(1-p)t}.$$
Taking $K$ as the support of $\mu$, we have
$$\int_0^{pt} v_{\delta t}(t-s,B_s){\rm d}A_s^{M\mu}
\leq e^{(-\Lambda_{\delta}+\varepsilon)(1-p)t}A_{pt}^{M\mu}.$$
Noting that $D_{pt}^t=C_{pt}^{\delta t,t}$, 
we get by \eqref{eq:dist},  
\begin{equation}\label{eq:est-ex}
\begin{split}
0\leq E_x\left[e^{A_{pt}^{(Q-1)\mu}}\right]-E_x\left[e^{D_{pt}^{t}}\right]
&\leq e^{(-\Lambda_{\delta}+\varepsilon)(1-p)t}E_x\left[e^{A_{pt}^{(Q-1)\mu}}A_{pt}^{M\mu}\right]\\
&=e^{-\lambda pt}e^{(-\Lambda_{\delta}+\varepsilon)(1-p)t}e^{\lambda pt}E_x\left[e^{A_{pt}^{(Q-1)\mu}}A_{pt}^{M\mu}\right].
\end{split}
\end{equation}

By the same argument  as for \eqref{eq:bdd}, 
there exist  $c_1>0$ and $c_2>0$ such that 
for any $\varepsilon_2>0$, 
$$e^{(\lambda-\varepsilon_2)pt}E_x\left[e^{A_{pt}^{(Q-1)\mu}}A_{pt}^{M\mu}\right]
\leq E_x\left[\sup_{0\leq s\leq pt}\left(e^{(\lambda-\varepsilon_2)s}e^{A_s^{(Q-1)\mu}}\right)A_{pt}^{M\mu}\right]
\leq c(\varepsilon_2)(c_1+c_2t).$$
Since $\Lambda_{\delta}>0$, 
there exists $\varepsilon_2>0$ for any $\varepsilon\in (0,\Lambda_{\delta})$ such that 
$$c:=(\Lambda_{\delta}-\varepsilon)(1-p)-\varepsilon_2 p>0.$$
Then the last term of \eqref{eq:est-ex} is less than 
$$
c(\varepsilon_2)e^{(-\lambda +\varepsilon_2)pt}e^{(-\Lambda_{\delta}+\varepsilon)(1-p)t}(c_1+c_2 t)
=c(\varepsilon_2)e^{-\lambda pt}e^{-ct}(c_1+c_2t),$$
that is,
$$0\leq e^{\lambda pt}\left(E_x\left[e^{A_{pt}^{(Q-1)\mu}}\right]-E_x\left[e^{D_{pt}^{t}}\right]\right)
\leq c_2(\varepsilon_2)e^{-ct}(c_1+c_2 t)\rightarrow 0\ (t\rightarrow\infty).$$
Hence by \eqref{eq-fkf-asymp},  
\begin{equation*}
\begin{split}
e^{\lambda pt}E_x\left[e^{D_{pt}^{t}}\right]
&=e^{\lambda pt}E_x\left[e^{A_{pt}^{(Q-1)\mu}}\right]
+e^{\lambda pt}\left(E_x\left[e^{D_{pt}^{t}}\right]-E_x\left[e^{A_{pt}^{(Q-1)\mu}}\right]\right)\\
&\rightarrow h(x)\int_{{\mathbb R}^d}h(y)\,{\rm d}y \ (t\rightarrow\infty).
\end{split}
\end{equation*}
This completes the proof. \qed
\medskip

We are now in a position to prove the lower bound of \eqref{eq:point-1}.  
For any $p\in (0,1)$, 
we have by the Markov property,
$$
E_x\left[e^{D_t^t};|B_t|\geq \delta t\right]
=E_x\left[e^{D_{pt}^t}
E_{B_{pt}}\left[e^{D_{(1-p)t}^t};|B_{(1-p)t}|\geq \delta t\right]\right].
$$
Since $D_{(1-p)t}^t\geq 0$ for any $t\geq 0$, 
the last term above is greater than
\begin{equation}\label{eq:side}
E_x\left[e^{D_{pt}^t}
P_{B_{pt}}\left(|B_{(1-p)t}|\geq \delta t\right)\right]
\geq E_x\left[e^{D_{pt}^t}
\right]P_0\left(|B_{(1-p)t}|\geq \delta t\right)
\end{equation}
by \cite[Appendix A]{S18}. 
Then by Lemma \ref{lem:fk-d}, we have as $t\rightarrow\infty$,
\begin{equation}\label{eq:optimal-0}
\begin{split}
&E_x\left[e^{D_{pt}^{t}}
\right]P_0\left(|B_{(1-p)t}|\geq \delta t\right)\\
&\sim \frac{\omega_d}{(2\pi)^{d/2}}\left(\frac{\delta^2 t}{1-p}\right)^{(d-2)/2}\exp\left(-\lambda pt-\frac{\delta^2 t}{2(1-p)}\right)
h(x) \int_{{\mathbb R}^d}h(y)\,{\rm d}y.
\end{split}
\end{equation}
If we let $p=1-\delta/\sqrt{-2\lambda}$,
then the last term of \eqref{eq:optimal-0} becomes
$$
\frac{\omega_d(\sqrt{-2\lambda}\delta)^{(d-2)/2}}{(2\pi)^{d/2}}t^{(d-2)/2}e^{(-\lambda-\sqrt{-2\lambda}\delta)t}
h(x) \int_{{\mathbb R}^d}h(y)\,{\rm d}y.
$$
We thus get the lower bound of \eqref{eq:point-1} by \eqref{eq:fk-ex}.

\section{Proof of Theorem \ref{thm-critical}}

Our proof of Theorem \ref{thm-critical} is 
a refinement of that of  Theorem \ref{thm:main}.
　
\subsection{Proof of (i)}
Let $\delta=\sqrt{-\lambda/2}$ 
and let $\{t_n\}$ be a positive increasing sequence such that $t_n\rightarrow\infty$ as $n\rightarrow\infty$. 
Let $G(t)$ be a positive function on $(0,\infty)$. 
For any $n\geq 1$ and $\varepsilon>0$, 
we have by the same way as in \eqref{eq:chebyshev} and \eqref{eq:max-markov},
\begin{equation}\label{eq:chebyshev-1}
{\mathbf P}_x\left(\max_{t_n\leq s\leq t_{n+1}}Z_s^{\delta t_n}\geq G(t_n)\right)
\leq E_x\left[e^{A_{t_n}^{(Q-1)\mu}}
{\mathbf E}_{B_{t_n}}\left[\max_{0\leq s\leq t_{n+1}-t_n}Z_s^{\delta t_n}\right]\right]/G(t_n).
\end{equation}

Let $a(t)$ be a nonnegative function on $(0,\infty)$ such that 
$a(t)=o(t) \ (t\rightarrow\infty)$ and   $R(t):=\delta t-a(t)$. 
For $s\leq t$, let ${\mathbf B}_s^{(t),k}$ be the position at time $s$ of the $k$th particle alive at time $t$.
Since  
$$\max_{0\leq s\leq t_{n+1}-t_n}Z_s^{\delta t_n}
\leq \sum_{k=1}^{Z_{t_{n+1}-t_n}}{\bf 1}_{\{\sup_{0\leq s\leq t_{n+1}-t_n}|{\mathbf B}_s^{(t_{n+1}-t_n),k}|\geq \delta t_n\}},$$
we have  by the same argument  as for \eqref{eq:est-max},
\begin{equation}\label{eq:est-max-1}
\begin{split}
&{\mathbf E}_x\left[\max_{t_n\leq s\leq t_{n+1}}Z_s^{\delta t_n}\right]
\leq E_x\left[e^{A_{t_n}^{(Q-1)\mu}}E_{B_{t_n}}\left[e^{A_{t_{n+1}-t_n}^{(Q-1)\mu}};\sup_{0\leq s\leq t_{n+1}-t_n}|B_s|\geq \delta t_n\right]\right]\\
&=E_x\left[e^{A_{t_n}^{(Q-1)\mu}}
E_{B_{t_n}}\left[e^{A_{t_{n+1}-t_n}^{(Q-1)\mu}};\sup_{0\leq s\leq t_{n+1}-t_n}|B_s|\geq \delta t_n\right]; 
|B_{t_n}|\geq R(t_n)\right]\\
&+E_x\left[e^{A_{t_n}^{(Q-1)\mu}}E_{B_{t_n}}\left[e^{A_{t_{n+1}-t_n}^{(Q-1)\mu}};\sup_{0\leq s\leq t_{n+1}-t_n}|B_s|\geq \delta t_n\right];  
|B_{t_n}|< R(t_n)\right]={\rm (I)}+{\rm (II)}.
\end{split}
\end{equation}

In what follows, we suppose that 
\begin{itemize}
\item $t_{n+1}-t_n\rightarrow 0$ as $n\rightarrow \infty$;
\item $a(t_n)^2/(t_{n+1}-t_n)\rightarrow\infty$ as $n\rightarrow\infty$.
\end{itemize}
Then by Remark \ref{rem:asymp-fk}, 
\begin{equation}\label{eq:(i)-1}
{\rm (I)}
\leq E_x\left[e^{A_{t_n}^{(Q-1)\mu}}; 
|B_{t_n}|\geq R(t_n)\right]\sup_{x\in {\mathbb R}^d}E_x\left[e^{A_{t_{n+1}-t_n}^{(Q-1)\mu}}\right]
\asymp e^{\sqrt{-2\lambda}a(t_n)}t_n^{(d-1)/2} \quad (n\rightarrow\infty).
\end{equation}
By the Cauchy-Schwarz inequality, 
we have for any $x\in {\mathbb R}^d$ with $|x|\leq R(t_n)$ and for any constants $p,q>1$ with $1/p+1/q=1$,
\begin{equation}\label{eq:cs}
\begin{split}
&E_x\left[e^{A_{t_{n+1}-t_n}^{(Q-1)\mu}};\sup_{0\leq s\leq t_{n+1}-t_n}|B_s|\geq \delta t_n\right]
\leq E_x\left[e^{A_{t_{n+1}-t_n}^{(Q-1)\mu}};
\sup_{0\leq s\leq t_{n+1}-t_n}|B_s-x|\geq a(t_n)\right]\\
&\leq E_x\left[e^{pA_{t_{n+1}-t_n}^{(Q-1)\mu}}\right]^{1/p}
P_x\left(
\sup_{0\leq s\leq t_{n+1}-t_n}|B_s-x|\geq a(t_n)\right)^{1/q}\\
&\leq cP_0\left(\sup_{0\leq s\leq t_{n+1}-t_n}|B_s|\geq a(t_n)\right)^{1/q}.
\end{split}
\end{equation}
If $r(t)$ is a positive function on $(0,\infty)$ such that $r(t)^2/t\rightarrow \infty$ as 
$t\rightarrow +0$, then by \cite[Corollary 3.4]{Se17} and the change of variables,
\begin{equation*}
\begin{split}
&P_0\left(\sup_{0\leq s\leq t}|B_s|\geq r(t)\right)
=P_0\left(\sup_{0\leq s\leq t/r(t)^2}|B_s|\geq 1\right)\\
&\asymp \int_0^{t/r(t)^2}\frac{e^{-1/(2t)}}{t^{(d+2)/2}}\,{\rm d}t
=\int_{r(t)^2/t}^{\infty}e^{-u/2}u^{(d-2)/2}\,{\rm d}u\sim 2e^{-r(t)^2/(2t)}\left(\frac{r(t)^2}{t}\right)^{(d-2)/2} 
\quad (t\rightarrow\infty).
\end{split}
\end{equation*}
Hence 
\begin{equation}\label{eq:bessel}
P_0\left(\sup_{0\leq s\leq t_{n+1}-t_n}|B_s|\geq a(t_n)\right)
\asymp \exp\left(-\frac{a(t_n)^2}{2(t_{n+1}-t_n)}\right) \left(\frac{a(t_n)^2}{t_{n+1}-t_n}\right)^{(d-2)/2} 
\quad (n\rightarrow\infty).
\end{equation}
For any $x\in {\mathbb R}^d$, since it follows by \eqref{eq-fkf-asymp} that 
$$E_x\left[e^{A_{t_n}^{(Q-1)\mu}}\right]\asymp e^{(-\lambda) t_n} \quad (n\rightarrow\infty),$$ 
we have by \eqref{eq:cs} and \eqref{eq:bessel},
\begin{equation}\label{eq:ii}
\begin{split}
{\rm (II)}
&\leq cE_x\left[e^{A_{t_n}^{(Q-1)\mu}}\right]P_0\left(\sup_{0\leq s\leq t_{n+1}-t_n}|B_s|\geq a(t_n)\right)^{1/q}\\
&\asymp e^{(-\lambda) t_n} \exp\left(-\frac{a(t_n)^2}{2q(t_{n+1}-t_n)}\right)\left(\frac{a(t_n)^2}{t_{n+1}-t_n}\right)^{(d-2)/(2q)}
 \quad (n\rightarrow\infty).
\end{split}
\end{equation}

For $c_1>0$, $c_2>0$ and $\alpha\in (0,1)$, if we let 
$$a(t)\equiv c_1, \quad t_n=c_2 n^{\alpha},$$
then
\begin{equation}\label{eq:comp-bdd}
c_2\alpha(n+1)^{\alpha-1}\leq t_{n+1}-t_n\leq c_2\alpha n^{\alpha-1}
\end{equation}
and therefore,
$$\frac{c_1}{\sqrt{c_2\alpha}}n^{(1-\alpha)/2}
\leq \frac{a(t_n)}{\sqrt{t_{n+1}-t_n}}\leq  \frac{c_1}{\sqrt{c_2\alpha}} (n+1)^{(1-\alpha)/2}$$
Since $t_{n+1}-t_n\rightarrow 0$ and $a(t_n)^2/(t_{n+1}-t_n)\rightarrow \infty$ as $n\rightarrow\infty$, 
we obtain by \eqref{eq:chebyshev-1}, \eqref{eq:est-max-1} \eqref{eq:(i)-1} and \eqref{eq:ii},
\begin{equation}\label{eq:chebyshev-1-2}
\begin{split}
&{\mathbf P}_x\left(\max_{t_n\leq s\leq t_{n+1}}Z_s^{\delta t_n}\geq G(t_n)\right)\\
&\leq \frac{c}{G(t_n)}
\left(e^{\sqrt{-2\lambda}a(t_n)}t_n^{(d-1)/2}+e^{(-\lambda) t_n} \exp\left(-\frac{a(t_n)^2}{2q(t_{n+1}-t_n)}\right)\left(\frac{a(t_n)^2}{t_{n+1}-t_n}\right)^{(d-2)/(2q)}\right)\\
&\leq \frac{c'}{G(t_n)}\left(e^{\sqrt{-2\lambda}c_1}n^{\alpha(d-1)/2}
+e^{c_2(-\lambda) n^{\alpha}}e^{-c_1^2n^{1-\alpha}/(2qc_2\alpha)}n^{(1-\alpha)(d-2)/(2q)}\right).
\end{split}
\end{equation}

Here we take $\alpha=1/2$ and $c_1$ so large that $c_1\geq \sqrt{-2q \lambda}c_2$.  
If we let $G(t)=t^a(\log t)(\log \log t)^{1+\varepsilon}$ for $a>0$ and $\varepsilon>0$, 
then by \eqref{eq:chebyshev-1-2}, 
$${\mathbf P}_x\left(\max_{t_n\leq s\leq t_{n+1}}Z_s^{\delta t_n}\geq G(t_n)\right)
\leq \frac{cn^{(d-1)/4}}{n^{a/2}(\log n)(\log \log n)^{1+\varepsilon}}.$$
In particular, if we let $a=(d+3)/2$, then  
$$
\sum_{n=1}^{\infty}{\mathbf P}_x\left(\max_{t_n\leq s\leq t_{n+1}}Z_s^{\delta t_n}\geq G(t_n)\right)<\infty.
$$
Hence by the Borel-Cantelli lemma, there exists an event of full ${\mathbf P}_x$-probability 
and a natural valued random variable $N\geq 1$ such that on this event,  
we have for all $n\geq N$,
$$\max_{t_n\leq s\leq t_{n+1}}Z_s^{\delta t_n}\leq G(t_n).$$
Moreover, for any $n\geq N$ and $t\in [t_n, t_{n+1}]$,
$$Z_t^{\delta t}\leq \max_{t_n\leq s\leq t_{n+1}}Z_s^{\delta t_n}\leq G(t_n)\leq G(t),$$
which completes the proof. 
\qed

\subsection{Proof of (ii)}

As in the proof of Lemma \ref{lem:lower-growth}, 
we denote by ${\mathbf B}_s^{t,k}$ the position of a particle at time $s$
starting from ${\mathbf B}_t^k$ at time $t$ ($s\geq t$). 
Let $\{t_n\}$ be a positive increasing sequence such that 
$t_n\rightarrow\infty$. 
Fix a constant $p_n\in [0,1)$ and a compact set $K\subset {\mathbb R}^d$. 
Then for each index $k$, 
\begin{equation*}
\begin{split}
&\left\{\text{${\mathbf B}_{p_nt_n}^{p_nt_n,k}\in K$, $|{\mathbf B}_s^{p_n t_n,k}|>\delta s$ for all $s\in [t_n,t_{n+1}]$}\right\}\\
&\supset 
\left\{\text{\begin{minipage}[c]{100mm}
${\mathbf B}_{p_n t_n}^{p_n t_n,k}\in K$, $|{\mathbf B}_{t_n}^{p_n t_n,k}|>|{\mathbf B}_{t_n}^{p_n t_n,k}-{\mathbf B}_{p_n t_n}^{p_n t_n,k}|
>\delta t_{n+1}+1$, 
$\sup_{t_n\leq s\leq t_{n+1}}|{\mathbf B}_s^{p_n t_n,k}-{\mathbf B}_{t_n}^{p_n t_n,k}|<1$
\end{minipage}}\right\}
=:E_n^k.
\end{split}
\end{equation*}
Let $G(t)$ and $f(t)$ be  positive functions on $(0,\infty)$ such that $f(t)\rightarrow 0$ as $t\rightarrow \infty$. 
Define  
$$N_t=\left\{e^{\lambda t}Z_t(K)\geq f(t)\right\}.$$
Then by the Markov property,
\begin{equation}\label{eq:c-markov-prop}
{\mathbf P}_x\left(\left\{\sum_{k=1}^{Z_{p_n t_n}}{\bf 1}_{E_n^k}\leq G(t_n)\right\}\cap N_{p_nt_n}\right)
={\mathbf E}_x\left[{\mathbf P}_{{\mathbf B}_{p_n t_n}}
\left(\sum_{k=1}^l {\bf 1}_{F_n^k}\leq G(t_n)\right)|_{l={Z_{p_n t_n}}}; 
N_{p_nt_n}\right]
\end{equation}
for 
$$F_n^k:=\left\{\text{\begin{minipage}[c]{100mm}
${\mathbf B}_0^{0,k}\in K$, 
$|{\mathbf B}_{(1-p_n)t_n}^{0,k}|>|{\mathbf B}_{(1-p_n)t_n}^{0,k}-{\mathbf B}_0^{0,k}|
>\delta t_{n+1}+1$, 
$\sup_{(1-p_n)t_n\leq s\leq t_{n+1}-p_n t_n}|{\mathbf B}_s^{0,k}-{\mathbf B}_{(1-p_n)t_n}^{0,k}|<1$
\end{minipage}}\right\}.$$ 
Let ${\mathbf x}=(x^1,\dots,x^l)$. Then by the same way as in \eqref{eq:indep-0} and \eqref{eq:indep},
\begin{equation}\label{eq:c-indep-0}
{\mathbf P}_{{\mathbf x}}
\left(\sum_{k=1}^l {\bf 1}_{F_n^k}\leq G(t_n)\right)
\leq e^{G(t_n)}\prod_{1\leq k\leq m, \, x^k\in K}\exp\left(-(1-e^{-1}){\mathbf P}_{x^k}(F_n^k)\right).
\end{equation}

Let
$$C_n:=\left\{|B_{(1-p_n)t_n}|>|B_{(1-p_n)t_n}-B_0|>\delta t_{n+1}+1\right\}.$$
Then for any $x\in K$, we have by the Markov property and the spatial uniformity of the Brownian motion,
\begin{equation*}
\begin{split}
{\mathbf P}_x(F_n^k)
&=P_x\left(\left\{\sup_{(1-p_n)t_n\leq s\leq t_{n+1}-p_nt_n}|B_s-B_{(1-p_n)t_n}|<1\right\}\cap C_n
\right)\\
&=E_x\left[ 
P_{B_{(1-p_n)t_n}}\left(\sup_{0\leq s\leq t_{n+1}-t_n}|B_s-B_0|<1\right); C_n\right]\\
&=P_x(C_n)
P_0\left(\sup_{0\leq s\leq t_{n+1}-t_n}|B_s-B_0|<1\right).
\end{split}
\end{equation*}

In what follows, we suppose that 
\begin{itemize}
\item $t_{n+1}-t_n\rightarrow 0$ as $n\rightarrow\infty$;
\item $t_{n+1}/\sqrt{(1-p_n)t_n}\rightarrow\infty$ as $n\rightarrow\infty$.
\end{itemize}
Then by the same way as in \eqref{eq:prob-c}, 
there exist $c_0>0$ and $c_1>0$ such that for any $x\in {\mathbb R}^d$, 
$$P_x(C_n)
\geq c_0\int_{(\delta t_{n+1}+1)/\sqrt{(1-p_n)t_n}}^{\infty}e^{-r^2/2}r^{d-1}\,{\rm d}r
\geq c_1 e^{-(\delta t_{n+1})^2/(2(1-p_n)t_n)}
\left(\frac{\delta t_{n+1}}{\sqrt{(1-p_n)t_n}}\right)^{d-2},$$
which implies that  for any $x\in K$,
$${\mathbf P}_x(F_n^k)
\geq c_1 e^{-(\delta t_{n+1})^2/(2(1-p_n)t_n)}\left(\frac{\delta t_{n+1}}{\sqrt{(1-p_n)t_n}}\right)^{d-2}
=:c_1 q_n.$$
Because of  this and \eqref{eq:c-indep-0}, there exists $c_2>0$ such that   
$$
{\mathbf P}_{{\mathbf x}}
\left(\sum_{k=1}^l {\bf 1}_{F_n^k}\leq G(t_n)\right)
\leq e^{G(t_n)}
\exp\left(-c_2q_n\cdot \sharp\left\{k\mid x^k\in K\right\}\right).
$$
Hence by \eqref{eq:c-markov-prop},
\begin{equation}\label{eq:upper-prob}
\begin{split}
{\mathbf P}_x\left(\left\{\sum_{k=1}^{Z_{p_n t_n}}{\bf 1}_{E_n^k}\leq G(t_n)\right\}
\cap N_{p_nt_n}\right)
&\leq e^{G(t_n)} 
{\mathbf E}_x\left[\exp\left(-c_2 q_n Z_{p_n t_n}(K)\right); 
N_{p_nt_n}\right]\\
&\leq \exp\left(G(t_n)-c_2 q_n e^{-\lambda p_n t_n}f(t_n)\right).
\end{split}
\end{equation}

Here we note that  
$$q_n e^{-\lambda p_n t_n}
=e^{g_n(p_n)}\left(\frac{\delta t_{n+1}}{\sqrt{(1-p_n)t_n}}\right)^{d-2}$$
for 
$$g_n(p)=-\lambda p t_n-\frac{(\delta t_{n+1})^2}{2(1-p)t_n} \quad (0\leq p<1).$$ 
Then the right hand side above takes the maximal value 
$g_n(p_n^*)=\lambda(t_{n+1}-t_n)$ for $p_n^*=1-t_{n+1}/(2t_n)$. 

We take $p_n=p_n^*$, 
$f(t)=(\log \log t)^{-\varepsilon}$ and 
$G(t)=c_3t^{(d-2)/2}(\log\log t)^{-\varepsilon}$ for $\varepsilon>0$ and $c_3>0$. 
Then 
\begin{equation}\label{eq:exponent}
G(t_n)-c_2 q_n e^{-\lambda p_n t_n}f(t_n)
=t_n^{(d-2)/2}(\log\log t_n)^{-\varepsilon}
\left\{c_3-c_2 e^{\lambda(t_{n+1}-t_n)}
\left(-\lambda t_{n+1}/t_n\right)^{(d-2)/2}\right\}.
\end{equation}
For some $\alpha\in (0,1]$ and $c>0$, if we let  
$$t_n=cn^{\alpha},$$
then $t_{n+1}-t_n\rightarrow 0$ as $n\rightarrow\infty$ 
and $t_{n+1}/\sqrt{(1-p_n)t_n}=\sqrt{2 t_{n+1}}\rightarrow\infty$ as $n\rightarrow\infty$. 
Therefore by \eqref{eq:exponent}, we can take $c_3>0$ so small that 
$$G(t_n)-c_2 q_n e^{-\lambda p_n t_n}f(t_n)
\leq -c_4n^{\alpha(d-2)/2}(\log\log n)^{-\varepsilon}$$
for some $c_4>0$. 
Then by \eqref{eq:upper-prob}, 
$$
{\mathbf P}_x\left(
\left\{\sum_{k=1}^{Z_{p_n t_n}}{\bf 1}_{E_n^k}\leq G(t_n)\right\}\cap N_{p_nt_n}\right)
\leq \exp\left(-c_4n^{\alpha(d-2)/2}(\log\log n)^{-\varepsilon}\right).
$$
In particular, if we assume that $d\geq 3$, then  
$$\sum_{n=1}^{\infty}
{\mathbf P}_x\left(
\left\{\sum_{k=1}^{Z_{p_n t_n}}{\bf 1}_{E_n^k}\leq G(t_n)\right\}\cap N_{p_nt_n}\right)<\infty.$$
Hence by the Borel-Cantelli lemma, the event 
$$\left\{\sum_{k=1}^{Z_{p_n t_n}}{\bf 1}_{E_n^k}>G(t_n)\right\}\cup (N_{p_nt_n})^c$$
occurs for all sufficiently large $n$. 
Since
$p_nt_n=t_n-t_{n+1}/2\rightarrow\infty$ as $n\rightarrow\infty$,
we have by \eqref{eq-limit-thm},
$$e^{\lambda p_n t_n} Z_{p_n t_n}(K)\rightarrow M_{\infty}\int_K h(y)\,{\rm d}y, \quad \text{${\mathbf P}_x$-a.s.}$$
Since $\int_K h(y)\,{\rm d}y>0$ and $f(t)\rightarrow 0$ as $t\rightarrow\infty$, 
we see that on the event   $\{M_{\infty}>0\}$, 
the event $(N_{p_nt_n})^c$ occurs only for finite $n\geq 1$, that is, 
the event $\left\{\sum_{k=1}^{Z_{p_n t_n}}{\bf 1}_{E_n^k}>G(t_n)\right\}$ occurs for all sufficiently large $n\geq 1$.

For all sufficiently large  $t>0$, there exists $n=n(t)\in {\mathbb N}$ such that $t_n\leq t<t_{n+1}$ 
and 
$$G(t_n)
=c_3t_n^{(d-2)/2}(\log\log t_n)^{-\varepsilon}
\geq c_5t_{n+1}^{(d-2)/2}(\log\log t)^{-\varepsilon}
\geq c_5t^{(d-2)/2}(\log\log t)^{-\varepsilon}$$
for some $c_5>0$. 
We thus have, ${\mathbf P}_x(\cdot\mid M_{\infty}>0)$-a.s.\ for all sufficiently large $t>0$, 
\begin{equation*}
\begin{split}
Z_t^{\delta t}
&=\sum_{k=1}^{Z_t}{\bf 1}_{\{|{\mathbf B}_t^k|>\delta t\}}
\geq \sum_{k=1}^{Z_{p_n t_n}}
{\bf 1}_
{\left\{\text{${\mathbf B}_{p_n t_n}^{p_n t_n,k}\in K$, 
$|{\mathbf B}_s^{p_n t_n,k}|>\delta s$ for all $s\in [t_n,t_{n+1}]$}
\right\}}\\
&\geq \sum_{k=1}^{Z_{p_n t_n}}{\bf 1}_{E_n^k}\geq G(t_n)
\geq c_5t^{(d-2)/2}(\log\log t)^{-\varepsilon}.
\end{split}
\end{equation*}
Since $\varepsilon>0$ is arbitrary, the last inequality above is valid 
by taking $c_5=1$. 
\qed

\appendix
\section{Appendix}
\subsection{Decay rate of the ground state}\label{appendix:decay}

Let  $\mu=\mu^+-\mu^-$ for some $\mu^+, \mu^- \in {\cal K}_{\infty}(1)$.   
Recall that  
$$\lambda(\mu)=\inf\left\{\frac{1}{2}\int_{{\mathbb R}^d}|\nabla u|^2\,{\rm d}x-\int_{{\mathbb R}^d}u^2\,{\rm d}\mu 
\mid u\in C_0^{\infty}({\mathbb R}^d), \int_{{\mathbb R}^d}u^2\,{\rm d}x=1\right\}.$$
In what follows, we let $\lambda:=\lambda(\mu)$ and assume that $\lambda<0$. 
As mentioned in Subsection \ref{subsection:fk}, 
$\lambda$ is the principal eigenvalue of the operator $-\Delta/2-\mu$ 
and the corresponding eigenfunction $h$ 
has a bounded, continuous and strictly positive version. 

\begin{lem}\label{lem:decay-ground}
Suppose that $\mu=\mu^+-\mu^-$ for some $\mu^+, \mu^- \in {\cal K}_{\infty}(1)$. 
Then for any positive constants $A_1$ and $A_2$ with $A_2<\sqrt{-2\lambda}<A_1$, 
there exist positive constants $c_1$, $c_2$ such that 
$$\frac{c_1e^{-A_1|x|}}{|x|^{(d-1)/2}}\leq h(x)
\leq \frac{c_2 e^{-A_2|x|}}{|x|^{(d-1)/2}}, \quad (|x|\geq 1).$$
Moreover, if $\mu^-$ is compactly supported in ${\mathbb R}^d$, 
then the inequality above holds with $A_1=\sqrt{-2\lambda}$. 
A similar result is valid for  $\mu^{+}$ and $A_2$. 
\end{lem}

\pf \ We follow the argument of \cite{C89} and \cite[Lemma 4.1]{T08}. 
We first discuss the upper bound of $h$. 
For $r>0$, let 
$$\sigma_r:=\inf\{t>0 \mid |B_t|\leq r\}.$$
Since $M_t:=e^{\lambda t}e^{A_t^{\mu}}h(B_t)$ is a $P_x$-martingale, 
we have by the optional stopping theorem and the H\"older inequality, 
\begin{equation}\label{eq:eigen-upper-1}
h(x)=E_x\left[e^{\lambda (t\wedge \sigma_r)}e^{A_{t\wedge \sigma_r}^{\mu}}h(B_{t\wedge \sigma_r})\right]
\leq \|h\|_{\infty}E_x\left[e^{p(\lambda+\varepsilon)(t\wedge \sigma_r)}\right]^{1/p}
E_x\left[e^{-q\varepsilon(t\wedge \sigma_r)}e^{qA_{t\wedge\sigma_r}^{\mu}}\right]^{1/q}
\end{equation}
for any $\varepsilon\in (0,-\lambda)$ and $p,q>1$ with $1/p+1/q=1$. 

Let $\mu_r({\rm d}x)={\bf 1}_{|x|>r}(x)\mu({\rm d}x)$ 
and let $\hat{P}_x$ be the law of the killed process of ${\mathbf M}$ 
by the exponential distribution with rate $q\varepsilon$. 
Then 
\begin{equation}\label{eq:eigen-upper-2}
E_x\left[e^{-q\varepsilon(t\wedge \sigma_r)}e^{qA_{t\wedge\sigma_r}^{\mu}}\right]
=E_x\left[e^{-q\varepsilon(t\wedge \sigma_r)}e^{qA_{t\wedge\sigma_r}^{\mu_r}}\right]
=\hat{E}_x\left[e^{qA_{t\wedge\sigma_r}^{\mu_r}}\right]
\leq\hat{E}_x\left[e^{qA_{\sigma_r}^{\mu_r^+}}\right].
\end{equation}
Since $\mu^+\in {\cal K}_{\infty}(1)$ and 
$$\hat{E}_x\left[qA_{\sigma_r}^{\mu_r^+}\right]
\leq q\hat{E}_x\left[A_{\infty}^{\mu_r^+}\right]
=q\int_{|y|\geq r}G_{q\varepsilon}(x,y)\,\mu^+({\rm }dy),$$
there exists $R=R(\varepsilon,p)>0$ such that for any $r\geq R$, 
$$\sup_{x\in {\mathbb R}^d}\hat{E}_x\left[qA_{\sigma_r}^{\mu_r^+}\right]
\leq q\sup_{x\in {\mathbb R}^d}\int_{|y|\geq r}G_{q\varepsilon}(x,y)\,\mu^+({\rm }dy)<1.$$
Then the Khasminskii lemma (see, e.g., \cite[Lemma 3.7]{CZ95}) implies that  for any $r\geq R$,
\begin{equation}\label{eq:exp-int}
\sup_{x\in {\mathbb R}^d}\hat{E}_x\left[e^{qA_{\sigma_r}^{\mu_r^+}}\right]<\infty.
\end{equation}
Hence by \eqref{eq:eigen-upper-1} and \eqref{eq:eigen-upper-2},
\begin{equation}\label{eq:upper-eigen-3}
\begin{split}
h(x)
&\leq \|h\|_{\infty}
E_x\left[e^{p(\lambda+\varepsilon)(t\wedge \sigma_r)}\right]^{1/p}\left(\sup_{x\in {\mathbb R}^d}
\hat{E}_x\left[e^{qA_{\sigma_r}^{\mu_r^+}}\right]\right)^{1/q}\\
&\rightarrow \|h\|_{\infty}
E_x\left[e^{p(\lambda+\varepsilon)\sigma_r}\right]^{1/p}
\left(\sup_{x\in {\mathbb R}^d}\hat{E}_x\left[e^{qA_{\sigma_r}^{\mu_r^+}}\right]\right)^{1/q} \ (t\rightarrow\infty).
\end{split}
\end{equation}

Let $\nu_r$ be the equilibrium potential of  
$B_r:=\{x\in {\mathbb R}^d \mid |x|\leq r\}$ (see \cite[p.82]{FOT11} for definition).
Then 
\begin{equation}\label{eq:potential}
E_x\left[e^{p(\lambda+\varepsilon)\sigma_r}\right]
=\int_{|y|\leq r}G_{-p(\lambda+\varepsilon)}(x,y)\nu_r({\rm d}y)
\leq \sup_{|y|\leq r}G_{-p(\lambda+\varepsilon)}(x,y)\nu_r(B_r).
\end{equation}
Since we see by \eqref{eq:resolvent} that 
for any $x,y\in {\mathbb R}^d$ with $|x|\geq 2r$ and $|y|\leq r$,
$$G_{-p(\lambda+\varepsilon)}(x,y)\leq c_{\varepsilon,p,r}\frac{e^{-\sqrt{-2p(\lambda+\varepsilon)}|x|}}{|x|^{(d-1)/2}},$$
we have by \eqref{eq:potential}, 
\begin{equation}\label{eq:potential-1}
E_x\left[e^{p(\lambda+\varepsilon)\sigma_r}\right]
\leq c_{\varepsilon,p,r}'\frac{e^{-\sqrt{-2p(\lambda+\varepsilon)}|x|}}{|x|^{(d-1)/2}} \ (|x|\geq 1).
\end{equation}
Then by \eqref{eq:upper-eigen-3},
\begin{equation}\label{eq:upper-eigen-4}
h(x)\leq \|h\|_{\infty}c_{\varepsilon,p,r}''
\left(\frac{e^{-\sqrt{-2p(\lambda+\varepsilon)}|x|}}{|x|^{(d-1)/2}}\right)^{1/p} \ 
(|x|\geq 1),
\end{equation}
which implies the desired upper bound of $h$. 
If we further assume that $\mu^+$ is compactly supported in ${\mathbb R}^d$, 
then  \eqref{eq:upper-eigen-4} is valid for $p=1$ 
because $\mu_r^+$ vanishes for large $r>0$.

We next discuss the lower bound of $h$. 
Here we denote by $\tilde{P}_x$ the law of the killed process of ${\mathbf M}$ 
by the exponential distribution with rate $-\lambda$. 
Then by the optional stopping theorem again,
\begin{equation}\label{eq:lower-eigen-1}
h(x)=E_x\left[e^{\lambda (t\wedge \sigma_r)}e^{A_{t\wedge \sigma_r}}h(B_{t\wedge \sigma_r})\right]
\geq \inf_{|y|\leq r}h(y)\tilde{E_x}\left[e^{-A_{t\wedge \sigma_r}^{\mu^-}}\right]
\geq \inf_{|y|\leq r}h(y)\tilde{E_x}\left[e^{-A_{\sigma_r}^{\mu^-}}; \sigma_r\leq t\right].
\end{equation}
By the same argument as for \eqref{eq:exp-int}, 
there exists $R=R(p)>0$ for any $p>1$ such that for any $r\geq R$, 
$$\sup_{x\in {\mathbb R}^d}{\tilde{E}_x\left[e^{A_{\sigma_r}^{\mu_r^{-}}/(p-1)}\right]^{p-1}}<\infty.$$
Then by the H\"older inequality, we have for any $p>1$ and $r\geq R$,
\begin{equation}\label{eq:lower-eigen-2}
\begin{split}
\tilde{E_x}\left[e^{-A_{\sigma_r}^{\mu^-}}; \sigma_r\leq t\right]
\geq \frac{\tilde{P}_x(\sigma_r\leq t)^p}{\tilde{E}_x\left[e^{A_{\sigma_r}^{\mu_r^{-}}/(p-1)}\right]^{p-1}}
&\geq \frac{\tilde{P}_x(\sigma_r\leq t)^p}{\sup_{y\in {\mathbb R}^d}\tilde{E}_y\left[e^{A_{\sigma_r}^{\mu_r^{-}}/(p-1)}\right]^{p-1}}\\
&\rightarrow \frac{\tilde{P}_x(\sigma_r<\infty)^p}
{\sup_{y\in {\mathbb R}^d}\tilde{E}_y\left[e^{A_{\sigma_r}^{\mu_r^{-}}/(p-1)}\right]^{p-1}} 
\quad (t\rightarrow\infty).
\end{split}
\end{equation}
Since
$$\tilde{P}_x(\sigma_r<\infty)=E_x\left[e^{\lambda \sigma_r};\sigma_r<\infty\right]=\int_{B_r}G_{-\lambda}(x,y)\nu_r({\rm d}y),$$ 
we have by the same argument as for \eqref{eq:potential-1}, 
$$\tilde{P}_x(\sigma_r<\infty)\geq c_{r}\frac{e^{-\sqrt{-2\lambda}|x|}}{|x|^{(d-1)/2}}.$$
Then by \eqref{eq:lower-eigen-1} and \eqref{eq:lower-eigen-2},
\begin{equation}\label{eq:lower-eigen}
h(x)\geq c_{p,r}
\left(\frac{e^{-\sqrt{-2\lambda}|x|}}{|x|^{(d-1)/2}}\right)^p \ 
(|x|\geq 1).
\end{equation}
We thus get the desired lower bound of $h$. 
If $\mu^-$ is compactly supported in ${\mathbb R}^d$, 
then \eqref{eq:lower-eigen} is valid for $p=1$
because $\mu_r^-$ vanishes for large $r>0$. 
\qed
\medskip

By using Lemma \ref{lem:decay-ground} instead of \cite[Lemma 4.1]{T08}, 
we can follow the argument of \cite[Section 4]{T08} to get  
\begin{thm}
Suppose that $\mu=\mu^+-\mu^-$ for some $\mu^+, \mu^- \in {\cal K}_{\infty}(1)$. 
Then for any $f\in {\cal B}_b({\mathbb R}^d)$,
$$\lim_{t\rightarrow\infty}e^{\lambda t}E_x\left[e^{A_t^{\mu}}f(B_t)\right]
=h(x)\int_{{\mathbb R}^d}f(y)h(y)\,{\rm d}y 
\quad (x\in {\mathbb R}^d).$$
\end{thm}

\subsection{Positivity of $M_{\infty}$ and survival}\label{appendix:positive}
We discuss relations among the positivity of $M_{\infty}$, 
the finiteness of the total number of branching and 
the survival property.  
Note that we already discussed in \cite[Proposition 3.6, Theorem 3.7 and Remark 3.14]{S08} 
the relation between the first and third properties 
for branching symmetric stable processes with absorbing boundary. 

Let $\overline{\mathbf M}=(\{{\mathbf B}_t\}_{t\geq 0}, \{\mathbf {P}_{\mathbf x}\}_{{\mathbf x}\in {\mathbf X}})$ 
be a branching Brownian motion on ${\mathbf X}$ 
with branching rate $\mu\in {\cal K}_{\infty}$ and branching mechanism $\{p_n(x)\}_{n\geq 0}$.  
Denote by $G^{\mu}(x,y)$ the Green function 
associated with the Feynman-Kac semigroup $p_t^{\mu}f(x)=E_x\left[e^{-A_t^{\mu}}f(B_t)\right]$.
For a function $u$ on ${\mathbb R}^d$, define 
$$F(u)(x)=\sum_{n=0}^{\infty}p_n(x)u(x)^n$$
if the right hand side makes sense. 
We first study the solution to the next equation:
\begin{equation}\label{eq:martingale}
u(x)=E_x\left[e^{-A_{\infty}^{\mu}}\right]
+E_x\left[\int_0^{\infty}e^{-A_t^{\mu}}F(u)(B_t)\,{\rm d}A_t^{\mu}\right], 
\quad 0\leq u(x)\leq 1 \quad (x\in {\mathbb R}^d).
\end{equation}

\begin{lem}\label{lem:equiv-1}
Suppose that 
\begin{equation}\label{eq:finite-energy}
\int_{{\mathbb R}^d}\int_{{\mathbb R}^d}G^{\mu}(x,y)\,\mu({\rm d}y)\mu({\rm d}x)<\infty.
\end{equation}
Let  $u$ and $v$ be functions on ${\mathbb R}^d$ 
such that $0\leq u(x)\leq v(x)<1$ on ${\mathbb R}^d$. 
If these functions are solutions to the equation \eqref{eq:martingale}, 
then $u\equiv v$.
\end{lem}

We omit the proof of Lemma \ref{lem:equiv-1} 
because it is similar to that of \cite[Lemma 3.5]{S08}. 
We note that if $\mu({\mathbb R}^d)<\infty$, then \eqref{eq:finite-energy} is fulfilled because 
$$\int_{{\mathbb R}^d}G^{\mu}(x,y)\,\mu({\rm d}y)=E_x\left[\int_0^{\infty}e^{-A_t^{\mu}}\,{\rm d}A_t^{\mu}\right]
=1-E_x\left[e^{-A_{\infty}^{\mu}}\right]\leq 1.$$

We next reveal the relations as we mentioned at the first of this subsection. 
Let $N$ be the total number of branching for $\overline{{\mathbf M}}$. 

\begin{prop}\label{prop:equiv}
Suppose that \eqref{eq:finite-energy} holds and ${\mathbf P}_x(M_\infty>0)>0$. 
If $d=1, 2$, then 
$$\{e_0=\infty\}=\{N=\infty\}=\{M_{\infty}>0\}, \quad \text{${\mathbf P}_x$-a.s.}$$
On the other hand, if $d\geq 3$, then
$$\{e_0=\infty\}\supsetneq\{N=\infty\}=\{M_{\infty}>0\}, \quad \text{${\mathbf P}_x$-a.s.}$$
\end{prop}

If $d\geq 3$, then the Brownian motion is transient 
so that the associated particle goes to infinity eventually. 
Since we assume that the branching rate $\mu$ is small at infinity, 
the number of branching can be small even on the survival event. 
In fact, branching never occurs with positive probability. 
On the other hand, if $d=1$ or $2$, then the Brownian motion is recurrent so that 
the associated particle can come to the support of $\mu$ infinitely often.   
Therefore, branching occurs infinite times on the survival event. 
\medskip

\pf \ Let $u(x)={\mathbf P}_x(N<\infty)$ and $v(x)={\mathbf P}_x(M_{\infty}=0)$.  
Then $v(x)<1$ by assumption. 
Moreover,  if $N<\infty$, then $Z_t$ is a finite random constant eventually 
and thus $M_{\infty}=0$. Namely, we obtain  $0\leq u(x)\leq v(x)<1$. 
Since $u$ and $v$ are solutions to the equation \eqref{eq:martingale}, 
we obtain $u\equiv v$ by Lemma \ref{lem:equiv-1}, whence
$$\{N=\infty\}=\{M_{\infty}>0\}, \quad \text{${\mathbf P}_x$-a.s.}$$

Let 
$u_e(x)={\mathbf P}_x(e_0<\infty)$. 
Then 
$$u_e(x)=E_x\left[\int_0^{\infty}e^{-A_t^{\mu}}
F(u_e)(B_t)\,{\rm d}A_t^{\mu}\right].$$
For $d=1,2$, since $P_x(A_{\infty}^{\mu}=\infty)=1$ by \cite[p.426, Proposition 3.11]{RY99}, 
we have $E_x[e^{-A_{\infty}^{\mu}}]=0$  
so that $u_e$ also satisfies the equation \eqref{eq:martingale}. 
Furthermore, since $\{e_0<\infty\}\subset\{M_{\infty}=0\}$, 
we obtain $0\leq u_e(x)\leq v(x)<1$ and thus  $u_e(x)=v(x)$ by Lemma \ref{lem:equiv-1}. 
This implies that 
$$\{e_0=\infty\}=\{M_{\infty}>0\}, \quad \text{${\mathbf P}_x$-a.s.}$$

On the other hand, if $d\geq 3$, 
then $\sup_{x\in {\mathbb R}^d}E_x[A_{\infty}^{\mu}]<\infty$ by \eqref{eq-green-bd}. 
Hence by Jensen's inequality, 
$${\mathbf P}_x(T=\infty)=E_x[e^{-A_{\infty}^{\mu}}]\geq \exp\left(-E_x[A_{\infty}^{\mu}]\right)
\geq \exp\left(-\sup_{x\in {\mathbb R}^d}E_x[A_{\infty}^{\mu}]\right)>0.$$
Since
$$\{e_0<\infty\}\cup \{T=\infty\}\subset \{N<\infty\}, \quad \{e_0<\infty\}\cap \{T=\infty\}=\emptyset,$$
we have $${\mathbf P}_x(N<\infty)\geq {\mathbf P}_x(e_0<\infty)+{\mathbf P}_x(T=\infty)
>{\mathbf P}_x(e_0<\infty).$$
Then by assumption,
$${\mathbf P}_x(e_0=\infty)>{\mathbf P}_x(N=\infty)={\mathbf P}_x(M_{\infty}>0)>0,$$
which shows that 
$$\{e_0=\infty\}\supsetneq \{N=\infty\}, \quad \text{${\mathbf P}_x$-a.s.}$$
We thus complete the proof. 
\qed

\subsection{Proof of \eqref{eq:int-upper-1}}\label{appendix:evaluate}
We  evaluate the integral in the right hand side of \eqref{eq:int-upper-1}. 
We first recall that $a(t)$ is a function on $(0,\infty)$ such that $a(t)=o(t) \ (t\rightarrow\infty)$ 
and $R(t)=\delta t+a(t)$ for some $\delta>0$. 
We will show that as $t\rightarrow\infty$,
\begin{equation}\label{eq:int-asymp}
\begin{split}
&(R(t)-\varepsilon_1 t)^d\int_0^t 
e^{(-\lambda+\varepsilon_2)s}\frac{1}{(t-s)^{(d+2)/2}}
\exp\left(-\frac{(R(t)-\varepsilon_1 t)^2}{2(t-s)}\right) \,{\rm d}s\\
&\asymp
\begin{cases}
e^{-(\delta-\varepsilon_1)^2t/2}t^{(d-2)/2} & (\delta>\sqrt{-2\lambda}),\\
e^{(-\lambda+\varepsilon_2)t-\sqrt{2(-\lambda+\varepsilon_2)}(R(t)-\varepsilon_1 t)}t^{(d-1)/2}
& (\delta\leq \sqrt{-2\lambda}).
\end{cases}
\end{split}
\end{equation}

By the change of variables $s=t-v$, we get 
\begin{equation}\label{eq:c-v}
\begin{split}
&\int_0^t e^{(-\lambda+\varepsilon_2)s}
\frac{1}{(t-s)^{(d+2)/2}}
\exp\left(-\frac{(R(t)-\varepsilon_1 t)^2}{2(t-s)}\right)\,{\rm d}s\\
&=e^{(-\lambda+\varepsilon_2)t}\int_0^t e^{-(-\lambda+\varepsilon_2)v}
\frac{e^{-(R(t)-\varepsilon_1 t)^2/(2v)}}{v^{(d+2)/2}}\,{\rm d}v\\
&=e^{(-\lambda+\varepsilon_2)t-\sqrt{2(-\lambda+\varepsilon_2)}(R(t)-\varepsilon_1 t)}
\int_0^t e^{-(\sqrt{(-\lambda+\varepsilon_2)v}-(R(t)-\varepsilon_1 t)/\sqrt{2v})^2}
\frac{1}{v^{(d+2)/2}}\,{\rm d}v\\
&=e^{(-\lambda+\varepsilon_2)t-\sqrt{2(-\lambda+\varepsilon_2)}(R(t)-\varepsilon_1 t)}{\rm (III)}.
\end{split}
\end{equation}
If we let 
$$w=\sqrt{(-\lambda+\varepsilon_2)v}-\frac{R(t)-\varepsilon_1 t}{\sqrt{2v}},$$
then  
\begin{equation}\label{eq:c-v-2}
({\rm III})=2\int_{-\infty}^{S(t)} e^{-w^2}
\frac{F_t(w)^d}{\sqrt{w^2+2\sqrt{2(-\lambda+\varepsilon_2)}(R(t)-\varepsilon_1 t)}}\,{\rm d}w
\end{equation}
for 
$$F_t(w)=\frac{2\sqrt{-\lambda+\varepsilon_2}}{w+\sqrt{w^2+2\sqrt{2(-\lambda+\varepsilon_2)}(R(t)-\varepsilon_1 t)}}$$
and
$$S(t)=\sqrt{(-\lambda+\varepsilon_2)t}-\frac{R(t)-\varepsilon_1 t}{\sqrt{2t}}.
$$

Assume that $\delta>\sqrt{-2\lambda}$. 
Fix $\varepsilon_1\in (0,\delta-\sqrt{-2\lambda})$, and 
take $\varepsilon_2>0$ so that $\delta-\varepsilon_1>\sqrt{2(-\lambda+\varepsilon_2)}$. 
Then there exist $c_1>0$ and $c_2>0$ such that for any $w\leq 0$ and $t>0$,
\begin{equation}\label{eq:asymp-1}
c_1\frac{\sqrt{w^2+2\sqrt{2(-\lambda+\varepsilon_2)}(R(t)-\varepsilon_1 t)}}{R(t)-\varepsilon_1 t} 
\leq F_t(w)
\leq c_2\frac{\sqrt{w^2+2\sqrt{2(-\lambda+\varepsilon_2)}(R(t)-\varepsilon_1 t)}}
{R(t)-\varepsilon_1 t}.
\end{equation}
Since $S(t)<0$ for all sufficiently large $t>0$, 
we have as $t\rightarrow\infty$,
\begin{equation*}
\begin{split}
({\rm III})
&\asymp \frac{1}{(R(t)-\varepsilon_1 t)^d}\int_{-\infty}^{S(t)}
e^{-w^2}\left(w^2+2\sqrt{2(-\lambda+\varepsilon_2)}(R(t)-\varepsilon_1 t)\right)^{(d-1)/2}\,{\rm d}w\\
&=\frac{1}{(R(t)-\varepsilon_1 t)^d}\int_{-S(t)}^{\infty} 
e^{-w^2}\left(w^2+2\sqrt{2(-\lambda+\varepsilon_2)}(R(t)-\varepsilon_1 t)\right)^{(d-1)/2}\,{\rm d}w.
\end{split}
\end{equation*}
If $w\geq -S(t)$, then 
$$w^2\geq \left(\frac{R(t)-\varepsilon_1 t}{\sqrt{2}t}-\sqrt{(-\lambda+\varepsilon_2)}\right)^2t
=\left(\frac{\delta-\varepsilon_1}{\sqrt{2}}-\sqrt{(-\lambda+\varepsilon_2)}+\frac{a(t)}{\sqrt{2}t}\right)^2t$$
and hence for all sufficiently large $t>0$,
$$c_1(\varepsilon_1)w^2\leq w^2+2\sqrt{2(-\lambda+\varepsilon_2)}(R(t)-\varepsilon_1 t)
\leq c_2(\varepsilon_1)w^2.$$
This implies that as $t\rightarrow\infty$,
\begin{equation}\label{eq:iii-1}
\begin{split}
({\rm III})
&\asymp \frac{c_3(\varepsilon_1)}{(R(t)-\varepsilon_1 t)^d}
\int_{-S(t)}^{\infty}e^{-w^2}w^{d-1}\,{\rm d}w
\sim \frac{c_3(\varepsilon_1)}{(R(t)-\varepsilon_1 t)^d}e^{-S(t)^2}(-S(t))^{d-2}\\
&=\frac{c_4(\varepsilon_1)}{(R(t)-\varepsilon_1 t)^d}
\exp\left(-\frac{(R(t)-\varepsilon_1 t)^2}{2t}\right)
e^{-((-\lambda+\varepsilon_2)t-\sqrt{2(-\lambda+\varepsilon_2)}(R(t)-\varepsilon_1 t))}
t^{(d-2)/2}.
\end{split}
\end{equation}

We next assume that $\delta\leq \sqrt{-2\lambda}$. 
Then 
\begin{equation*}
\begin{split}
{\rm (III)}
&=\int_{-\infty}^0 e^{-w^2}
\frac{F_t(w)^d}{\sqrt{w^2+2\sqrt{2(-\lambda+\varepsilon_2)}(R(t)-\varepsilon_1 t)}}\,{\rm d}w\\
&+\int_0^{S(t)} e^{-w^2}
\frac{F_t(w)^d}{\sqrt{w^2+2\sqrt{2(-\lambda+\varepsilon_2)}(R(t)-\varepsilon_1 t)}}\,{\rm d}w
=({\rm III})_1+({\rm III})_2.
\end{split}
\end{equation*}
By \eqref{eq:asymp-1} and the change of variables ($w=-v$), 
we obtain as $t\rightarrow\infty$,
\begin{equation*}
\begin{split}
({\rm III})_1
&\asymp \frac{1}{(R(t)-\varepsilon_1 t)^d}\int_0^{\infty} 
e^{-v^2}\left(v^2+2\sqrt{2(-\lambda+\varepsilon_2)}(R(t)-\varepsilon_1 t)\right)^{(d-1)/2}\,{\rm d}v\\
&\asymp \frac{c_4(\varepsilon_1)}{(R(t)-\varepsilon_1 t)^{(d+1)/2}}.
\end{split}
\end{equation*}
Since there exist $c(\varepsilon_1)>0$ and $c'(\varepsilon_1)>0$ such that 
for all sufficiently large $t>0$,
$$c(\varepsilon_1)\sqrt{t}\leq 
w+\sqrt{w^2+2\sqrt{2(-\lambda+\varepsilon_2)}(R(t)-\varepsilon_1 t)}
\leq c'(\varepsilon_1)\sqrt{t} \quad (0\leq w\leq S(t))$$
and $S(t)\rightarrow\infty$ as $t\rightarrow\infty$,
we also have 
$$
({\rm III})_2\asymp \frac{c_5(\varepsilon_1)}{(R(t)-\varepsilon_1 t)^{(d+1)/2}} \quad (t\rightarrow\infty),$$
that is, 
\begin{equation}\label{eq:iii-2}
({\rm III})\asymp \frac{c_6(\varepsilon_1)}{(R(t)-\varepsilon_1 t)^{(d+1)/2}} \quad (t\rightarrow\infty).
\end{equation}
We thus get \eqref{eq:int-asymp} by \eqref{eq:c-v}, \eqref{eq:c-v-2}, \eqref{eq:iii-1} and \eqref{eq:iii-2}. 
\medskip

\noindent
{\bf Acknowledgments} \ 
The author would like to thank 
Professor Naotaka Kajino, Professor Tomoyuki Shirai, 
Professor Ryokichi Tanaka, and one of the referees of \cite{S18}
for their valuable comments motivating this work. 
He is grateful to Professor Xia Chen for his comment 
which improved Theorem \ref{thm-critical}.

\address{
Yuichi Shiozawa\\
Department of Mathematics\\
Graduate School of Science\\ 
Osaka University\\
Toyonaka, Osaka, 560-0043,
Japan
}
{\texttt{shiozawa@math.sci.osaka-u.ac.jp}}
\end{document}